\documentclass[10pt]{article}
\usepackage[utf8]{inputenc}
\usepackage[T1]{fontenc}
\usepackage[english]{babel}
\usepackage{indentfirst}
\usepackage{latexsym}
\usepackage[lmargin=2cm,tmargin=2.5cm,bmargin=2.5cm,rmargin=2.5cm]%
{geometry}
\usepackage{amsbsy, amsmath,amsfonts,amssymb,amsthm}
\usepackage{times}
\usepackage{graphicx}
\usepackage{amsmath}
\usepackage{amsfonts}
\usepackage{amssymb}%
\newtheorem{theorem}{Theorem}[section]
\newtheorem{definition}[theorem]{Definition}

\newtheorem{corollary}[theorem]{Corollary}
\newtheorem{lemma}[theorem]{Lemma}
\newtheorem{remark}[theorem]{Remark}

\numberwithin{equation}{section}
\usepackage[pagewise]{lineno}

\begin{document}

\title{On the fractional heat semigroup and product estimates in Besov spaces and applications in theoretical analysis of the fractional Keller-Segel system}

\author{\\{{Jhean E. P\'erez-L\'opez}{\thanks{Email: jelepere@uis.edu.co (corresponding author)}},\hspace{0.3cm}
\ {Diego A. Rueda-G\'omez} {\thanks{Email: diaruego@uis.edu.co}},\hspace{0.3cm}
\ {\'Elder J. Villamizar-Roa} {\thanks{Email: jvillami@uis.edu.co}}}\\\\{\small Universidad Industrial de Santander} \\{\small {Escuela de Matem\'aticas, A.A. 678, Bucaramanga, Colombia}}}
\date{}
\maketitle

\begin{abstract}
This paper is concerned with the fractional Keller-Segel system in the temporal and spatial variables. We consider fractional dissipation for the physical variables including a fractional dissipation mechanism for the chemotactic diffusion, as well as a time fractional variation assumed in the Caputo sense. We analyze the fractional heat semigroup obtaining time decay and integral estimates of the Mittag-Leffler operators in critical Besov spaces, and prove a bilinear estimate derived from the nonlinearity of the Keller-Segel system, without using auxiliary norms. We use these results in order to prove the existence of global solutions in critical homogeneous Besov spaces employing only the norm of the natural persistence space, including the existence of self-similar solutions, which constitutes a persistence result in this framework. In addition, we prove a uniqueness result without assuming any smallness condition of the initial data.

\bigskip{} \noindent\textbf{Keywords:} Fractional Keller-Segel system, Besov spaces,  bilinear estimate,  global existence, uniqueness.

\medskip{} \noindent\textbf{AMS MSC:} 35R11, 35A01, 35K55, 35Q92, 92C17.

\end{abstract}

\section{Introduction}\label{intro}
We consider the Keller-Segel system which  describes the movement of living organisms towards higher concentration regions of chemical attractants.  This system is composed of two coupled parabolic equations describing the interaction between the density of cells and the concentration of the chemoattractant, which reads as follows:
\begin{equation}
	\left\{ \begin{array}{ll}
		\eta_{t}-D_{\eta}\Delta\eta=-\chi\nabla\cdot(\eta\nabla v), & \mbox{ in }\mathbb{R}^{n}\times(0,\infty),\\
		v_{t}-D_{v}\Delta v=-\gamma v+\kappa\eta, & \mbox{ in }\mathbb{R}^{n}\times(0,\infty),\\
		\eta(x,0)=\eta_{0}(x),\ v(x,0)=v_{0}(x), & \mbox{ in }\mathbb{R}^{n},
	\end{array}\right.\label{eq: Keller-Segel intro_cl}
\end{equation}
where $n\geq 1.$ In (\ref{eq: Keller-Segel intro_cl}), $\eta, v$ are the unknowns denoting the density of cells and the chemical concentration, respectively. The parameters $D_\eta$ and $D_v$ represent the corresponding diffusion coefficients for $\eta$ and $v,$ while $\chi,$ $\gamma$ and $\kappa$ are nonnegative parameters denoting the chemotactic sensitivity, and the decay and production rates, respectively. The issues of existence and long-time asymptotic behaviour of solutions for the Keller-Segel system have attracted the attention of many authors (see for instance \cite{Corrias,Duarte,Ferreira,Hillen,Nagai,Winkler1,Winkler2} and references therein). In particular, for $n=1$ it is known the existence of global solutions, and the blow-up is entirely ruled out. In two and three dimensions, there exist global solutions for small data. In the two-dimensional case, solutions of (\ref{eq: Keller-Segel intro_cl}) with total mass of class $m<4\pi$ remain bounded for all times, while for $\epsilon>0,$ there exist unbounded solutions with total mass of cells $m<4\pi+\epsilon.$ For $n\geq 3,$ system (\ref{eq: Keller-Segel intro_cl}) has unbounded solutions for arbitrarily small mass of cells \cite{Winkler2}.\\

The classical Keller-Segel model assumes that the density diffusion is not affected by the nonlocal behaviour of the organisms. However, in many situations found in nature, organisms develop alternative search strategies, particularly when chemoattractants, food, or other targets are sparse or rare.  Then, the trajectories of the population of organisms are better described by the so called L\'evy flights than Brownian motion (see \cite{Escudero,Klafter}). L\'evy flights  behaviour has been suggested
in numerous biological contexts, including  immune cells, ecology, and human
populations (c.f. \cite{Estrada} and references for a deeper discussion). This consideration motivates the substitution of the classical diffusion in the Keller-Segel system (\ref{eq: Keller-Segel intro_cl}) by a fractional diffusion. On the other hand, regarding to the flux by chemotaxis, it is also relevant to consider that the attraction force be replaced by a less singular interaction kernel. This last consideration has been point out relevant in the analysis of the propagation of chaos for some aggregation-diffusion models \cite{Salem}. In addition, taking into account that the behavior of most biological systems has memory properties, which are neglected when an integer-order time derivative is assumed, we also assume a time variation in a fractional framework. This introduces a
nonlocal delay in time for the moving population \cite{Estrada2}.
Based on observations such as those mentioned, we are interested in the theoretical analysis of the following Keller-Segel system in the fractional setting
\begin{equation}
	\left\{ \begin{array}{ll}
		^cD^\alpha_{t}\eta+D_{\eta}(-\Delta)^{\theta/2}\eta=-\chi\nabla\cdot(\eta G(v)), & \mbox{ in }\mathbb{R}^{n}\times(0,\infty),\\
		^cD^\alpha_{t}v+D_{v}(-\Delta)^{\theta/2}v=-\gamma v+\kappa\eta, & \mbox{ in }\mathbb{R}^{n}\times(0,\infty),\\
		\eta(x,0)=\eta_{0}(x),\ v(x,0)=v_{0}(x), & \mbox{ in }\mathbb{R}^{n},
	\end{array}\right.\label{eq: Keller-Segel intro}
\end{equation}
where $^cD^\alpha_{t}$ denotes the time fractional derivative operator of order $\alpha\in(0,1)$ in the Caputo sense.   

We recall that if $f\in L^1(0,T;X),$ $T>0,$ and $X$ is a Banach space, the Riemann-Liouville fractional integral of order $\alpha$ of $f$ is defined by 
$$I^\alpha_tf(t)=\frac{1}{\Gamma(\alpha)}\int_0^t(t-\tau)^{\alpha-1}f(\tau)d\tau,\ t\in [0,T]. $$
In addition, if $f\in C([0,T];X), 0<T\leq \infty,$ is such that $I^{1-\alpha}_tf\in W^{1,1}(0,T;X)$, the Caputo fractional derivative of order $\alpha$ of $f$ is defined by
$$^cD^\alpha_{t}f(t):=\frac{d}{dt}\left\{  I^{1-\alpha}_t[f(t)-f(0)]\right\}=\frac{d}{dt}\left\{ \int_0^t(t-\tau)^{-\alpha}[f(\tau)-f(0)]d\tau\right\}.$$
In (\ref{eq: Keller-Segel intro}), $(-\Delta)^{\theta/2},$ $\theta\in (0,2],$ denotes the fractional laplacian operator of order $\theta/2$ defined  by $(-\Delta)^{\theta/2}f(x)=\mathcal{F}^{-1}(\vert \xi\vert^\theta \hat{f}(\xi))(x),$ where $\hat{f}(\xi)=\mathcal{F}(f)(\xi)$ and $\mathcal{F}^{-1}(f)(\xi)$ denote the Fourier transform and the inverse Fourier transform of $f$, respectively. In addition, $G(v)$ is also a nonlocal term defined by 
\[
\ensuremath{G(v)(x)=\nabla\left((-\Delta)^{-\theta_{1}/2}v\right)(x),\quad x\in\mathbb{R}^{n},}
\]
for $\theta_{1}\in[0,n),$ which can be alternatively represented by $G(v)f=K(x)\ast f,$ $K(x)\sim \frac{x}{\vert x\vert^{n-\theta_1}}.$  The case $\theta=2$ and $\theta_1=0$ corresponds to the classical Keller-Segel system (\ref{eq: Keller-Segel intro_cl}). For $\theta_1=0,\alpha=1$ and $n=2,$ in \cite{Biler} the authors proved a result of local existence and uniqueness of solution for (\ref{eq: Keller-Segel intro}) in homogeneous Besov spaces by using some estimates of the linear dissipative equation in the framework of mixed temporal-spatial spaces, the Chemin mono-norm methods, the Fourier localization  and the Littlewood-Paley theory. Later, in \cite{Zhai}, the author proved the existence, uniqueness and stability of solutions for (\ref{eq: Keller-Segel intro}) in critical Besov spaces under smallness condition on the initial data. The results of \cite{Zhai}, are based on the $L^p$-$L^q$ time decay for the semigroup $e^{-t(-\Delta)^{\theta/2}}$ in Besov spaces, which leads to use auxiliar norms of Besov type and Kato-time-weighted norms. Some results of global existence and blow-up for the particular case of (\ref{eq: Keller-Segel intro}) with  $\alpha=1$ and without considering the $v$-equation, have been obtained in  \cite{Bournaveas,Shi} and some references therein.\\

On the other hand, the time fractional Keller-Segel system has not been extensively studied. In \cite{Cuevas1} the authors studied the global existence and long time behaviour of solutions for the particular case of (\ref{eq: Keller-Segel intro}) assuming $\theta=2, \theta_1=0$ and $\alpha\in(0,1),$ with small initial data in the Besov-Morrey space $\mathcal{N}^{-b}_{r,\lambda,\infty}\times \dot{B}_{\infty,\infty}$ with $n\geq 2,$ $0\leq \lambda \leq n-2,$ $b=2-\frac{n-\lambda}{r}$ and $\frac{n-\lambda}{2}<r<n-\lambda,$ in the same spirit of the results of \cite{Ferreira} for the classical Keller-Segel system (\ref{eq: Keller-Segel intro_cl}).  Some regularity properties of solution for (\ref{eq: Keller-Segel intro}) assuming $\theta=2, \theta_1=0$ and $\alpha\in(0,1),$ with initial data in $L^n\cap L^{n/2}\cap L^\infty\times \dot{B}_{\infty,\infty}$ where recently obtained in \cite{Cuevas2}. It is worthwhile to observe that the existence space in \cite{Cuevas1} includes auxiliary norms {\it a la} Kato, like in \cite{Zhai}. In that approach, the fixed point argument is applied by considering a suitable time-dependent $X$ whose norm is given by the sum of a norm $L^\infty(0,\infty; X_1)$ and a norm of kind $\sup_{t>0}t^a\Vert u\Vert_{X_2},$ for some $a\neq 0.$ With this type of norm is possible to deal with the bilinear term coming from the cross-difussion term; however, the uniqueness in the natural class $C([0,T);\tilde{X}),$ where $\tilde{X}$ corresponds to the maximal closed subspace of $X$ in which the heat semigroup is continuous, is not possible.\\

Motivated by the above considerations, the aim of this paper is to analyze the existence, uniqueness and persistence of global solutions for the spatio-temporal fractional Keller-Segel system (\ref{eq: Keller-Segel intro}) in the framework of critical Besov spaces without using auxiliary norms. In order to get this aim, we first derive time decay and integral estimates of the Mittag-Leffler operators in critical Besov spaces, and prove a bilinear estimate derived from the nonlinearity of the Keller-Segel system employing only the norm of the natural persistence space. In order to estimate the bilinear operator, in addition to dealing with the action of the fractional heat semigroup, is necessary to prove a  product estimate in the homogeneous Besov setting.\\

In order to establish the main result, we start by recalling the mild formulation of (\ref{eq: Keller-Segel intro}) in the fractional setting. According to Duhamel's principle, the system (\ref{eq: Keller-Segel intro})
is formally equivalent to the following integral formulation:
\begin{equation} 
	\begin{cases}
		{ \eta(t)=E_{\alpha}(-t^{\alpha}(-\Delta)^{\theta/2})\eta_{0}-\int_{0}^{t}(t-\tau)^{\alpha-1}\nabla\cdot E_{\alpha,\alpha}(-(t-\tau)^{\alpha}(-\Delta)^{\theta/2})(\eta G(v))(\tau)d\tau,}\\
		{v(t)=E_{\alpha}(-t^{\alpha}((-\Delta)^{\theta/2}-\gamma))v_{0}+ \int_{0}^{t}(t-\tau)^{\alpha-1}E_{\alpha,\alpha}(-(t-\tau)^{\alpha}((-\Delta)^{\theta/2}-\gamma))\eta(\tau)d\tau.}
	\end{cases}\label{eq:Mild formulation}
\end{equation}

Here $\{E_\alpha(-t^\alpha(-\Delta)^{\theta/2})\}_{t \geq 0}$ and  $\{E_{\alpha,\alpha}(-t^\alpha(-\Delta)^{\theta/2})\}_{t \geq 0}$ denote the Mittag-Leffler families defined by
\begin{align*}
	& E_{\alpha}(-t^{\alpha}(-\Delta)^{\theta/2})=\int_{0}^{\infty}M_{\alpha}(\tau)U_{\theta}(\tau t^{\alpha})d\tau,\\
	& E_{\alpha,\alpha}(-t^{\alpha}(-\Delta)^{\theta/2})=\int_{0}^{\infty}\alpha \tau M_{\alpha}(\tau)U_{\theta}(\tau t^{\alpha})d\tau,
\end{align*}
where $U_{\theta}(t)$ is the fractional heat semigroup defined in Fourier variables as
$\widehat{U_{\theta}(t)f}=\ensuremath{e^{-t\left|\xi\right|^{\theta}}\widehat{f}},$ and $M_\alpha:\mathbb{C}\rightarrow \mathbb{C}$ is the Mainardi function which is defined by
$$M_\alpha(z)=\sum_{n=0}^\infty\frac{z^n}{n!\Gamma(1-\alpha(1+n))}.$$

In the classical case $\alpha=1,$ according to Duhamel's principle, the system (\ref{eq: Keller-Segel intro_cl})
is formally equivalent to the following integral formulation:
\begin{equation}
	\begin{cases}
		\eta(t) & =U_{\theta}(t)\eta_{0}-\displaystyle\int_{0}^{t}\nabla\cdot U_{\theta}(t-\tau)(\eta G(v))(\tau)d\tau,\\
		v(t) & =\widetilde{U}_{\theta}(t)v_{0}+\displaystyle\int_{0}^{t}\widetilde{U}_{\theta}(t-\tau)\eta(\tau)d\tau,
	\end{cases}\label{eq:Mild formulation_b}
\end{equation}
where $\widetilde{U}_{\theta}(t)=e^{-\gamma t}U_{\theta}(t).$ 

A solution
$\left[\eta,v\right]$ of the integral system (\ref{eq:Mild formulation})
is called a mild solution of the differential system (\ref{eq: Keller-Segel intro}). In the rest of this work, we will denote the
bilinear and linear operators appearing in (\ref{eq:Mild formulation}) as:
\begin{eqnarray}
	B_{\theta}(\eta,v)(t) & = & -\int_{0}^{t}(t-\tau)^{\alpha-1}\nabla\cdot E_{\alpha,\alpha}(-(t-\tau)^{\alpha}(-\Delta)^{\theta/2})(\eta G(v))(\tau)d\tau,\label{bil1}\\
	T_{\theta}(\eta)(t) & = & \int_{0}^{t}(t-\tau)^{\alpha-1}E_{\alpha,\alpha}(-(t-\tau)^{\alpha}((-\Delta)^{\theta/2}-\gamma))\eta(\tau)d\tau.\label{eq:Operadores B y T}
\end{eqnarray}
Thus, we rewrite system (\ref{eq:Mild formulation}) as follows
\begin{equation}
	\begin{cases}
		\eta(t) & =E_{\alpha}(-t^{\alpha}(-\Delta)^{\theta/2})\eta_{0}+B_{\theta}(\eta,v)(t),\\
		v(t) & =E_{\alpha}(-t^{\alpha}((-\Delta)^{\theta/2}-\gamma))v_{0}+T_{\theta}(\eta)(t).
	\end{cases}\label{eq: Mild formulation OPERADORES}
\end{equation}

Note that if $\gamma=0,$ the system (\ref{eq: Keller-Segel intro})
has a scaling property. Indeed, it is not difficult to check that if $[\eta,v]$
is a regular solution of $(\ref{eq: Keller-Segel intro})$ (with $\gamma=0$), then the pair $[\eta_{\sigma},v_{\sigma}]$
defined by 
\begin{equation}
	\eta_{\sigma}(x,t):=\sigma^{2\theta+\theta_{1}-2}\eta\left(\sigma x,\sigma^{\frac{\theta}{\alpha}}t\right)\mbox{ and }v_{\sigma}(x,t):=\sigma^{\theta+\theta_{1}-2}v\left(\sigma x,\sigma^{\frac{\theta}{\alpha}}t\right),\label{eq: Scaling}
\end{equation}
is also a solution of (\ref{eq: Keller-Segel intro}). In this case,
the map 
\begin{eqnarray}\label{scaling}
	[\eta,v]\longmapsto[\eta_{\sigma},v_{\sigma}],
\end{eqnarray}
is called the scaling of (\ref{eq: Keller-Segel intro}), and solutions
invariant by the scaling, this is, solutions $[\eta,v]$ such that $[\eta,v]=[\eta_{\sigma},v_{\sigma}]$ for all $\sigma>0$,
are called self-similar solutions. Note that if $[\eta,v]$ is a self-similar
solution, the initial data $[\eta_{0},v_{0}]$ must be invariant by
the scaling 
\begin{equation}
	[\eta_{0},v_{0}]\longmapsto [\eta_{\sigma0},v_{\sigma0}],\label{eq:Scaling dato inicial}
\end{equation}
and from (\ref{eq: Scaling}) it must have 
\[
\eta_{0}(x)=\sigma^{2\theta+\theta_{1}-2}\eta_{0}\left(\sigma x\right)\qquad\mbox{and}\qquad v_{0}(x)=\sigma^{\theta+\theta_{1}-2}v_{0}\left(\sigma x\right),
\]
this is, a necessary condition to obtain self-similar solutions is
that the data $\eta_{0}$ and $v_{0}$ be homogeneous functions
of degrees $2-2\theta-\theta_{1}$ and $2-\theta-\theta_{1},$ respectively.

In the case $\gamma\ne0$ the system (\ref{eq: Keller-Segel intro})
has not a scaling property; however, we can use the ``intrinsic scaling'' (\ref{scaling})
in order to choose the function spaces of initial data. Explicitly, we consider the following class of initial data (see the notations in Section 2):
\[
\eta_{0}\in\dot{B}{}_{p,\infty}^{2-2\theta-\theta_{1}+\frac{n}{p}}\qquad\mbox{and}\qquad v_{0}\in\dot{B}{}_{q,\infty}^{2-\theta-\theta_{1}+\frac{n}{q}}.
\]

Now, we are in position to establish the main results of this paper.
\begin{theorem} \label{lem: produto} (Product estimate) Let $n\ge1$, $\theta_{1}\in[0,n)$,
	$\frac{6n}{5n+\theta_{1}}<p\le q\le p^{\prime}$, $\max\left\{ 1,\ensuremath{1-\frac{n}{2}-\frac{\theta_{1}}{2}+\frac{n}{p}}\right\}$ $<\theta<1+\frac{n-\theta_{1}}{3}$,
	and $\rho_{1},\rho_{2}\ge0$ small enough. Then, for $f\in\dot{B}{}_{p,\infty}^{2-2\theta-\theta_{1}+\frac{n}{p}+\rho_{1}}$
	and $g\in\dot{B}{}_{q,\infty}^{2-\theta-\theta_{1}+\frac{n}{q}+\rho_{2}}$,
	we have that $fG(g)\in\dot{B}{}_{p,\infty}^{3-3\theta-\theta_{1}+\frac{n}{p}+\rho_{1}+\rho_{2}}$
	and 
	\begin{equation}
		\left\Vert fG(g)\right\Vert _{\dot{B}{}_{p,\infty}^{3-3\theta-\theta_{1}+\frac{n}{p}+\rho_{1}+\rho_{2}}}\leq C\left\Vert f\right\Vert _{\dot{B}{}_{p,\infty}^{2-2\theta-\theta_{1}+\frac{n}{p}+\rho_{1}}}\left\Vert g\right\Vert _{\dot{B}{}_{q,\infty}^{2-\theta-\theta_{1}+\frac{n}{q}+\rho_{2}}}.\label{eq:EstimativaProducto}
	\end{equation}	
\end{theorem}

\begin{theorem}{\label{Bilinearfinal}} (Bilinear estimate)
	Let $n\ge1$, $0<T\le\infty,$ $\theta_{1}\in[0,n)$, $\alpha>0$,
	$\frac{6n}{5n+\theta_{1}}<p\le q\le p^{\prime}$ and $\max\left\{ 1,\ensuremath{1-\frac{n}{2}-\frac{\theta_{1}}{2}+\frac{n}{p}}\right\} <\theta<1+\frac{n-\theta_{1}}{3}$.
	Then, there exists a constant $K>0$ (independent of $T$) such that
	\begin{equation}
		\|B_{\theta}(\eta,v)\|_{L^{\infty}\left((0,T);\dot{B}{}_{p,\infty}^{2-2\theta-\theta_{1}+\frac{n}{p}}\right)}\leq K\|\eta\|_{\left(L^{\infty}(0,T);\dot{B}{}_{p,\infty}^{2-2\theta-\theta_{1}+\frac{n}{p}}\right)}\|v\|_{L^{\infty}\left((0,T);\dot{B}{}_{q,\infty}^{2-\theta-\theta_{1}+\frac{n}{q}}\right)},\label{eq:EcuacionBilinear}
	\end{equation}
	for all $\eta\in L^{\infty}\left((0,T);\dot{B}{}_{p,\infty}^{2-2\theta-\theta_{1}+\frac{n}{p}}\right)$
	and $v\in L^{\infty}\left((0,T);\dot{B}{}_{q,\infty}^{2-\theta-\theta_{1}+\frac{n}{q}}\right).$
\end{theorem}

\begin{theorem} \label{Teo: Main} (Well-posedness)
	Let $n\ge 1$, $\theta_{1}\in[0,n)$, $\alpha>0$,
	$\frac{6n}{5n+\theta_{1}}<p\le q\le p^{\prime}$ and let $\theta$ such that  $\max\left\{ 1,\ensuremath{1-\frac{n}{2}-\frac{\theta_{1}}{2}+\frac{n}{p}}\right\}$ $<\theta<1+\frac{n-\theta_{1}}{3}$.
	There exist $\varepsilon>0$ and $\delta>0$ such that, if 
	\[
	\left\Vert \eta_{0}\right\Vert _{\dot{B}{}_{p,\infty}^{2-2\theta-\theta_{1}+\frac{n}{p}}}<\varepsilon\qquad \mbox{and}\qquad\left\Vert v_{0}\right\Vert _{\dot{B}{}_{q,\infty}^{2-\theta-\theta_{1}+\frac{n}{q}}}<\varepsilon,
	\]
	then there exists a unique mild solution $[\eta,v]$ for $(\ref{eq: Keller-Segel intro})$
	such that 
	\[
	\left\Vert \eta\right\Vert _{L^{\infty}\left((0,\infty);\dot{B}{}_{p,\infty}^{2-2\theta-\theta_{1}+\frac{n}{p}}\right)}<\delta\qquad \mbox{and}\qquad\left\Vert v\right\Vert _{L^{\infty}\left((0,\infty);\dot{B}{}_{q,\infty}^{2-\theta-\theta_{1}+\frac{n}{q}}\right)}<\delta.
	\]	
\end{theorem}

\begin{corollary} \label{Cor: Corolario 1} (Self-similarity)
	Assume the hypotheses of Theorem \ref{Teo: Main} with $\gamma=0,$ and consider $\eta_{0}\in\dot{B}{}_{p,\infty}^{2-2\theta-\theta_{1}+\frac{n}{p}}$
	and $v_{0}\in\dot{B}{}_{q,\infty}^{2-\theta-\theta_{1}+\frac{n}{q}}$
	being homogeneous functions with degrees $2-2\theta-\theta_{1}$ and
	$2-\theta-\theta_{1}$, respectively. Then the solution $[\eta,v]$ obtained
	through Theorem \ref{Teo: Main} is self-similar.
\end{corollary}
\begin{remark}
	\begin{enumerate}
		\item Theorem \ref{Bilinearfinal} plays a central role in the persistence part of Theorem \ref{Teo: Main} and is also central to the proof of Theorem \ref{Cor: Unicidadsinnorma} below. Moreover, to our knowledge, this type of bilinear estimate for the Keller-Segel system is new in the context of critical spaces.
		
		\item Theorem \ref{Teo: Main}, additionally to existence and uniqueness, establishes a persistence result because we do not use auxiliary norms in the solution spaces as it is used in previous works (cf. \cite{Zhai,Cuevas1}). In particular, considering $\theta\neq 2$ and $\theta_1\neq0,$ Theorem \ref{Teo: Main} complements the existence and uniqueness result in \cite{Cuevas1}, as well as for $\alpha\neq 1,$ Theorem \ref{Teo: Main} complements the existence and uniqueness result in \cite{Zhai}.
		
		\item The proof of Theorem \ref{Teo: Main} in the case $\alpha=1$ is carried out taking into account the mild formulation (\ref{eq:Mild formulation_b}). If we denote by $[\eta_\alpha,v_\alpha],$ the mild solution of (\ref{eq: Keller-Segel intro}) in the sense of (\ref{eq:Mild formulation}) for $\alpha\in(0,1),$ and $[\eta_1,v_1]$ the mild solution of (\ref{eq: Keller-Segel intro}) in the sense of (\ref{eq:Mild formulation_b}) for $\alpha=1,$ it is not clear if  $\lim_{\alpha\rightarrow 1^-}[\eta_\alpha,v_\alpha]=[\eta_1,v_1].$ This is an open question that beyond being raised in this model, can be formulated in general parabolic problems (cf. \cite{Planas}).
		
		\item The analysis carried out in the proof of Theorem \ref{Teo: Main} allows us to include negative values for the parameter $\theta_1,$ namely, $\theta_1\in (-2n,0)$ (cf. Lemma \ref{lem: produto} and Remark \ref{rem:aux1}). However, we do not know the possible physical meaning in the description of the model. 
		
	\end{enumerate}
\end{remark}

Using the estimates developed in the proof of Theorem
\ref{Teo: Main}, we prove the following uniqueness theorem without assuming any smallness condition of the initial data. The existence of solutions for arbitrary large initial data is an open problem. This uniqueness result seems new for chemotaxis problems in the context of critical spaces, including the classical Keller-Segel system (\ref{eq: Keller-Segel intro_cl}). This issue has been raised in the context of Navier-Stokes equations (see \cite{Jhean2} and some references therein).
\vspace{0.1cm}

In general, given a Banach space $X$ we denote by $\tilde{X}$ the
maximal closed subspace of $X$ in which the family of operators $\{E_{\alpha}(-t^{\alpha}(-\Delta)^{\theta/2})\}_{t\geq0}$
is continuous.

\begin{theorem}[Uniqueness] \label{Cor: Unicidadsinnorma} Let $n\ge1$,
	$0<T\le\infty,$ $\theta_{1}\in[0,n)$, $\alpha>0$, $\frac{6n}{5n+\theta_{1}}<p\le q\le p^{\prime}$
	and $\max\left\{ 1,\ensuremath{1-\frac{n}{2}-\frac{\theta_{1}}{2}+\frac{n}{p}}\right\} <\theta<1+\frac{n-\theta_{1}}{3}$.
	If $\left[\eta^{1},v^{1}\right]$ and $\left[\eta^{2},v^{2}\right]$
	are two mild solutions of (\ref{eq: Keller-Segel intro}) in $C\left([0,T);\dot{B}{}_{p,\infty}^{2-2\theta-\theta_{1}+\frac{n}{p}}\right)\times C\left([0,T);\dot{B}{}_{q,\infty}^{2-\theta-\theta_{1}+\frac{n}{q}}\right)$
	with the same initial data $\left[\eta_{0},v_{0}\right]\in\tilde{\dot{B}}{}_{p,\infty}^{2-2\theta-\theta_{1}+\frac{n}{p}}\times\tilde{\dot{B}}{}_{q,\infty}^{2-\theta-\theta_{1}+\frac{n}{q}}$,
	then $\left[\eta^{1}(t),v^{1}(t)\right]=\left[\eta^{2}(t),v^{2}(t)\right]$
	in $\dot{B}{}_{p,\infty}^{2-2\theta-\theta_{1}+\frac{n}{p}}\times\dot{B}{}_{q,\infty}^{2-\theta-\theta_{1}+\frac{n}{q}}$
	for all $t\in[0,T).$ 
\end{theorem}

The rest of this paper is organized as follows. In Section \ref{preliminaries}, we give some preliminaries about
Besov spaces. Section \ref{Keyestimates} is devoted to
the proof of the linear and nonlinear estimates; in particular, we prove the product and bilinear estimates established. Finally, in Section
\ref{proofs}, we prove our results about existence and uniqueness of mild solutions (\ref{eq:Mild formulation}).

\section{Preliminaries}\label{preliminaries}

Briefly we recall some preliminaries about Besov spaces. In what follows $\varphi$ denotes a radially symmetric function such
that 
\[
\varphi\in C_{c}^{\infty}\left(\mathbb{R}^{n}\backslash\left\{ 0\right\} \right),\,\mbox{supp}\,\varphi\subset\left\{ x\,;\,\frac{3}{4}\leq\left\vert x\right\vert \leq\frac{8}{3}\right\} ,
\]
and
\[
\sum\limits _{j\in\mathbb{Z}}\varphi_{j}(\xi)=1,\,\text{\ }\forall\xi\in\mathbb{R}^{n}\backslash\left\{ 0\right\} ,\text{ where }\varphi_{j}(\xi):=\varphi\left(\xi2^{-j}\right).
\]
Recall the localization operators $\Delta_{j}$ and $S_{k}$ defined
by
\[
\Delta_{j}f=\varphi_{j}(D)f=(\varphi_{j})^{\vee}\ast f\text{\ \,and\ \,}S_{k}f=\sum\limits _{j\leq k}\Delta_{j}f.\text{ }
\]
One can check easily the identities
\[
\Delta_{j}\Delta_{k}f=0\,\,\mbox{if}\,\,\left\vert j-k\right\vert \geq2\text{ and }\Delta_{j}\left(S_{k-2}g\Delta_{k}f\right)=0\,\,\mbox{if}\,\,\left\vert j-k\right\vert \geq5.
\]
Moreover, we have the Bony's decomposition (see \cite{Bony})
\begin{equation}
	fg=T_{f}g+T_{g}f+R(fg),\label{eq: descomposion de Bony}
\end{equation}
where
\[
T_{f}g=\sum\limits _{j\in\mathbb{Z}}S_{j-2}f\Delta_{j}g,\text{ }R(fg)=\sum\limits _{j\in\mathbb{Z}}\Delta_{j}f\tilde{\Delta}_{j}g\text{ }\ \text{and\ }\tilde{\Delta}_{j}g=\sum\limits _{\left\vert j-j^{\prime}\right\vert \leq1}\Delta_{j^{\prime}}g.
\]
We also denote $\tilde{\varphi}_{j}=\varphi_{j-1}+\varphi_{j}+\varphi_{j+1}$
and $\tilde{D}_{j}=D_{j-1}\cup D_{j}\cup D_{j+1}$ where $j\in\mathbb{Z}$
and $D_{j}=\left\{ x\,;\,\frac{3}{4}2^{j}\leq\left\vert x\right\vert \leq\frac{8}{3}2^{j}\right\} $.
Notice that $\tilde{\varphi}_{j}=1$ on $D_{j}$.\\

\begin{lemma} \label{lem:Bernstein inequality} {(}\textbf{{Bernstein
			inequality)} } \cite{Jhean} Assume that $1\le q\le p\le\infty$. Then
	\begin{equation}
		\left\Vert f\right\Vert _{L^{p}\left(\mathbb{R}^{n}\right)}\leq C2^{j\left(\frac{n}{q}-\frac{n}{p}\right)}\left\Vert f\right\Vert _{L^{q}\left(\mathbb{R}^{n}\right)},\label{eq:Bernstein Lebesgue}
	\end{equation}
	for all $f\in L^{q}\left(\mathbb{R}^{n}\right)$ such that $\mbox{supp}\hat{f}\subset D_{j}$.
	
\end{lemma}

\begin{definition} \label{Def: Definicion espacios Besov} Let $1\leq p,r\leq\infty$
	and $s\in\mathbb{R}$. The homogeneous Besov space $\dot{B}_{p,r}^{s}=\dot{B}_{p,r}^{s}\left(\mathbb{R}^{n}\right)$
	is defined as
	\[
	\dot{B}_{p,r}^{s}=\left\{ f\in\mathcal{S}^{\prime}(\mathbb{R}^{n})/\mathcal{P};\,\left\Vert f\right\Vert _{\dot{B}_{p,r}^{s}}<\infty\right\} ,
	\]
	where
	\begin{equation}
		\left\Vert f\right\Vert _{\dot{B}_{p,r}^{s}}:=\left\{ \begin{array}{l}
			\left(\sum\limits _{j\in\mathbb{Z}}2^{jsr}\left\Vert \Delta_{j}f\right\Vert _{L^{p}}^{r}\right)^{\frac{1}{r}}\,\,\mbox{if}\,\,r<\infty,\\
			\mathop{\sup}\limits _{j\in\mathbb{Z}}\,2^{js}\left\Vert \Delta_{j}f\right\Vert _{L^{p}}\,\,\,\,\,\,\,\,\,\,\,\,\,\,\,\,\mbox{if}\,\,r=\infty.
		\end{array}\right.\label{Norma Besov}
	\end{equation}
	
\end{definition}

\begin{lemma} \label{Lem: Interpolacion entre Besov para dar Besov.}\cite{Berg}
	Let $1\le p\le\infty,$ $1\le r,r_{0},r_{1}\leq\infty$ and $s,s_{0},s_{1}\in\mathbb{R}$
	be such that $s=\left(1-\theta\right)s_{0}+\theta s_{1}$ with $\theta\in\left(0,1\right)$.
	Then 
	\[
	\left(\dot{B}{}_{p,r_{0}}^{s_{0}},\dot{B}{}_{p,r_{1}}^{s_{1}}\right)_{\theta,r}=\dot{B}{}_{p,r}^{s}.
	\]
	
\end{lemma}

\begin{lemma} \label{Lem:Dualidad entre besov-morrey e besov--bloques}\cite{Berg}
	Let $1<p\le\infty,$ $1<r\leq\infty$ and $s\in\mathbb{R}$. Then
	\[
	\left(\dot{B}{}_{p^{\prime},r^{\prime}}^{-s}\right)^{\prime}=\dot{B}{}_{p,r}^{s}.
	\]
\end{lemma}

\section{Key estimates}\label{Keyestimates}

\subsection{Time decay estimates of Mittag-Leffler operators}
In order to estimate the bilinear operator $B_\theta \left(\cdot,\cdot\right)$
and the linear operator $T_\theta(\cdot)$ introduced in (\ref{bil1}) and (\ref{eq:Operadores B y T}),
we need to deal with the action of the operators $E_{\alpha}(-t^{\alpha}\left(-\Delta\right)^{\theta/2})$ and  $E_{\alpha,\alpha}(-t^{\alpha}\left(-\Delta\right)^{\theta/2})$
in Besov spaces.

\begin{lemma} \label{Lem: Heat-Kernel-estimate en Besov}\cite{Zhai} Let $\theta>0$ and $\zeta\ge0$ and consider the fractional heat operator $U_{\theta}(t)$ defined in Fourier variables as
	$\widehat{U_{\theta}(t)f}=\ensuremath{e^{-t\left|\xi\right|^{\theta}}\widehat{f}}.$ If $s_{1}\leq s_{2},$
	$1\leq p_{1}\leq p_{2}\leq\infty$ and $1\leq r\leq\infty$, then
	the following inequality holds
	\begin{eqnarray}
		\left\Vert (-\Delta)^{\zeta/2}U_{\theta}(t)f\right\Vert _{\dot{B}_{p_{2},r}^{s_{2}}}\leq Ct^{-\frac{s_{2}-s_{1}+\zeta}{\theta}-\frac{1}{\theta}\left(\frac{n}{p_{1}}-\frac{n}{p_{2}}\right)}\|f\|_{\dot{B}_{p_{1},r}^{s_{1}}}.\label{es4}
	\end{eqnarray}
	
\end{lemma}
\begin{lemma}\label{Minardi}\cite{Planas}
	Let $\alpha\in(0,1)$ and $-1<r<\infty.$ Then $M_\alpha(t)\geq 0$ for all $t\geq 0$ and
	$$\int_0^\infty t^rM_\alpha(t)dt=\frac{\Gamma(r+1)}{\Gamma(\alpha r+1)}.$$
\end{lemma}

\begin{lemma} \label{Lem: Frac_Heat-Kernel-estimate en Besov} Let
	$\theta>0$ and $\zeta\ge0$. 
	
	If $s_{1}\leq s_{2},$ $1\leq p_{1}\leq p_{2}\leq\infty,$ $1\leq r\leq\infty$
	and $\frac{1}{\theta}(s_{1}-s_{2}+\zeta+\frac{n}{p_{1}}-\frac{n}{p_{2}})<1,$
	then the following inequality holds 
	\begin{eqnarray}
		\left\Vert (-\Delta)^{\zeta/2}E_{\alpha}(-t^{\alpha}(-\Delta)^{\theta/2})f\right\Vert _{\dot{B}_{p_{2},r}^{s_{2}}}\leq Ct^{-\frac{\alpha}{\theta}(s_{2}-s_{1}+\zeta)-\frac{\alpha}{\theta}\left(\frac{n}{p_{1}}-\frac{n}{p_{2}}\right)}\|f\|_{\dot{B}_{p_{1},r}^{s_{1}}}.\label{es2}
	\end{eqnarray}
	Moreover, if $s_{1}\leq s_{2},$ $1\leq p_{1}\leq p_{2}\leq\infty,$
	$1\leq r\leq\infty$ and $\frac{1}{\theta}(s_{1}-s_{2}+\zeta+\frac{n}{p_{1}}-\frac{n}{p_{2}})<2,$
	then the following inequality holds 
	
	\begin{eqnarray}
		\left\Vert (-\Delta)^{\zeta/2}E_{\alpha,\alpha}(-t^{\alpha}(-\Delta)^{\theta/2})f\right\Vert _{\dot{B}_{p_{2},r}^{s_{2}}}\leq Ct^{-\frac{\alpha}{\theta}(s_{2}-s_{1}+\zeta)-\frac{\alpha}{\theta}\left(\frac{n}{p_{1}}-\frac{n}{p_{2}}\right)}\|f\|_{\dot{B}_{p_{1},r}^{s_{1}}}.\label{es2-1}
	\end{eqnarray}
	
\end{lemma} 

\begin{proof} Let $f\in\dot{B}_{p_{1},r}^{s_{1}}.$ From Lemma \ref{Minardi}
	and estimate (\ref{es4}) in Lemma \ref{Lem: Heat-Kernel-estimate en Besov},
	it holds 
	\begin{eqnarray*}
		&  & \left\Vert (-\Delta)^{\zeta/2}E_{\alpha}(-t^{\alpha}(-\Delta)^{\theta/2})f\right\Vert _{\dot{B}_{p_{2},r}^{s_{2}}}\leq\int_{0}^{\infty}M_{\alpha}(\tau)\left\Vert (-\Delta)^{\zeta/2}U_{\theta}(\tau t^{\alpha})f\right\Vert _{\dot{B}_{p_{2},r}^{s_{2}}}d\tau\\
		&  & \ \ \leq C\left[\int_{0}^{\infty}M_{\alpha}(\tau)\tau^{-\frac{1}{\theta}(s_{2}-s_{1}+\zeta)-\frac{1}{\theta}\left(\frac{n}{p_{1}}-\frac{n}{p_{2}}\right)}d\tau\right]t^{-\frac{\alpha}{\theta}(s_{2}-s_{1}+\zeta)-\frac{\alpha}{\theta}\left(\frac{n}{p_{1}}-\frac{n}{p_{2}}\right)}\|f\|_{\dot{B}_{p_{1},r}^{s_{1}}}\\
		&  & \ \ \leq Ct^{-\frac{\alpha}{\theta}(s_{2}-s_{1}+\zeta)-\frac{\alpha}{\theta}\left(\frac{n}{p_{1}}-\frac{n}{p_{2}}\right)}\|f\|_{\dot{B}_{p_{1},r}^{s_{1}}},\ t>0.
	\end{eqnarray*}
	Wich prove (\ref{es2}). The proof of (\ref{es2-1}) follows analogously. 
	
\end{proof}

\subsection{Integral Estimates}

In order to estimate the integral terms in (\ref{eq:Mild formulation})
we present a version in Besov spaces of the Yamazaki estimate obtained
in \cite{Yamazaki} in the context of Lorentz spaces $L^{(p,d)}$. In the contex of Besov-Lorentz-Morrey spaces and working in the non frational case, a related estimate was proved in \cite{Jhean2}. We remark that the estimate presented here is more general and we do not need (although it is possible) to use Lorentz spaces as base space for Besov, and therefore
our way of prove is different to that presented in \cite{Jhean2}.

\begin{lemma} \label{Lem: Estimativa de integral en Besov Para operador Bilineal}
	Let $1\le p\le\infty$, $\zeta\ge0$, $\alpha>0$, $\theta>\zeta$
	and $s_{0},s\in\mathbb{R}$ be such that $-s+\theta-\zeta=-s_{0}.$
	Then, there exists a constant $C>0$ such that 
	\[
	\int_{0}^{\infty}\left\Vert \tau^{\alpha-1}(-\Delta)^{\zeta/2}E_{\alpha,\alpha}(-\tau^{\alpha}(-\Delta)^{\theta/2})f\right\Vert _{\dot{B}{}_{p,1}^{-s_{0}}}d\tau\leq C\left\Vert f\right\Vert _{\dot{B}{}_{p,1}^{-s}},
	\]
	for all $f\in\dot{B}{}_{p,1}^{-s}.$
	
\end{lemma}

\begin{proof} Let $f\in\dot{B}{}_{p,1}^{-s}$ and define the function
	$h_{f}$ by 
	\[
	h_{f}\left(\tau\right)=\left\Vert \tau^{\alpha-1}(-\Delta)^{\zeta/2}E_{\alpha,\alpha}(-\tau^{\alpha}(-\Delta)^{\theta/2})f\right\Vert _{\dot{B}{}_{p,1}^{-s_{0}}}.
	\]
	Thus, for $-s_{i}\le-s_{0}$ $(i=1,2)$ and using Lemma \ref{Lem: Frac_Heat-Kernel-estimate en Besov}
	we have
	
	\begin{align*}
		h_{f}\left(\tau\right) & =\left\Vert \tau^{\alpha-1}(-\Delta)^{\zeta/2}E_{\alpha,\alpha}(-\tau^{\alpha}(-\Delta)^{\theta/2})f\right\Vert _{\dot{B}{}_{p,1}^{-s_{0}}}\\
		& \le C\tau^{-\frac{\alpha}{\theta}\left(-s_{0}-(-s_{i})+\zeta+\frac{(1-\alpha)}{\alpha}\theta\right)}\left\Vert f\right\Vert _{\dot{B}{}_{p,1}^{-s_{i}}.}
	\end{align*}
	Taking, for example, $-s_{1}=-s-\varepsilon$ and $-s_{2}=-s+\varepsilon$
	for $\varepsilon$ small enough, we have that $-s_{i}<-s_{0}$, and defining $\frac{1}{z_{i}}=\frac{\alpha}{\theta}\left(-s_{0}+s_{i}+\zeta\right)$
	we obtain 
	
	\begin{align*}
		\frac{1}{z_{1}} & =\frac{\alpha(-s_{0}+s_{1}+\zeta+\frac{(1-\alpha)}{\alpha}\theta)}{\theta}=\frac{\alpha(-s_{0}+s+\text{\ensuremath{\varepsilon}}+\zeta+\frac{(1-\alpha)}{\alpha}\theta)}{\theta}\\
		& =\frac{\alpha}{\theta}\left(\theta+\text{\ensuremath{\varepsilon}}+\frac{(1-\alpha)}{\alpha}\theta\right)=1+\frac{\alpha}{\theta}\text{\ensuremath{\varepsilon}}>1,
	\end{align*}
	and 
	\begin{align*}
		0<\frac{1}{z_{2}} & =\frac{\alpha(-s_{0}+s_{2}+\zeta+\frac{(1-\alpha)}{\alpha}\theta)}{\theta}=\frac{\alpha(-s_{0}+s-\text{\ensuremath{\varepsilon}}+\zeta+\frac{(1-\alpha)}{\alpha}\theta)}{\theta}\\
		& =\frac{\alpha}{\theta}\left(\theta-\text{\ensuremath{\varepsilon}}+\frac{(1-\alpha)}{\alpha}\theta\right)=1-\frac{\alpha}{\theta}\text{\ensuremath{\varepsilon}}<1.
	\end{align*}
	
	Therefore $0<z_{1}<1<z_{2}<\infty,$ and for $\phi=1/2$ we have $1=\frac{\phi}{z_{1}}+\frac{1-\phi}{z_{2}}$
	and $-s=(1-\phi)(-s_{1})+\phi(-s_{2}).$ Thus, for $i=1,2,$ it follows that $h_{f}\in L^{z_{i},\infty}\left(0,\infty\right)$
	with the estimate $\left\Vert h_{f}\right\Vert _{L^{z_{i},\infty}\left(0,\infty\right)}\leq C\left\Vert f\right\Vert _{\dot{B}{}_{p,1}^{-s_{i}}}$ and we can use interpolation in Lorentz spaces and Lemma \ref{Lem: Interpolacion entre Besov para dar Besov.}
	in order to get 
	\[
	\left\Vert h_{f}\right\Vert _{L^{1}\left(0,\infty\right)}\leq C\left\Vert f\right\Vert _{\dot{B}{}_{p,1}^{-s}},
	\]
	which finishes the proof.
\end{proof}

\begin{remark}
	
	Under the same conditions in Lemma \ref{Lem: Estimativa de integral en Besov Para operador Bilineal}
	we also have 
	\[
	\int_{0}^{\infty}\left\Vert \tau^{\alpha-1}(-\Delta)^{\zeta/2}E_{\alpha,\alpha}(-t^{\alpha}((-\Delta)^{\theta/2}-\gamma))f\right\Vert _{\dot{B}{}_{p,1}^{-s_{0}}}d\tau\leq C\left\Vert f\right\Vert _{\dot{B}{}_{p,1}^{-s}},
	\]
	for all $f\in\dot{B}{}_{p,1}^{-s}.$
	
\end{remark}

The next lemma concerns with an estimate for the operator $\mathcal{B}$
defined by 
\begin{eqnarray}
	\mathcal{B}\left(f\right):=\int_{0}^{\infty}\tau^{\alpha-1}\nabla\cdot E_{\alpha,\alpha}(-\tau^{\alpha}(-\Delta)^{\theta/2})f(\cdot,\tau)d\tau.\label{opB}
\end{eqnarray}

\begin{lemma} \label{Lem: Operador integral B} Let $1<p\le\infty$, $\alpha>0$, $\theta>1$
	and $s_{0},s\in\mathbb{R}$ be such that $-s+\theta-1=-s_{0}.$
	Then, there exists $C>0$ such that
	\[
	\left\Vert \mathcal{B}\left(f\right)\right\Vert _{\dot{B}{}_{p,\infty}^{s}}\leq C\mathop{\sup}\limits _{t>0}\left\Vert f(t)\right\Vert _{\dot{B}{}_{p,\infty}^{s_{0}}},
	\]
	for all $f\in L^{\infty}\left(\left(0,\infty\right);\dot{B}{}_{p,\infty}^{s_{0}}\right)$.
	
\end{lemma}

\begin{proof} Using duality and Lemma \ref{Lem: Estimativa de integral en Besov Para operador Bilineal} with $\zeta=1$
	we have 
	\begin{align}
		& \left\Vert \mathcal{B}\left(f\right)\right\Vert _{\dot{B}{}_{p,\infty}^{s}}=\mathop{\sup}\limits _{\left\Vert h\right\Vert _{\dot{B}{}_{p^{\prime},1}^{-s}}=1}\left\vert \left\langle \mathcal{B}\left(f\right),h\right\rangle \right\vert \leq C\mathop{\sup}\limits _{\left\Vert h\right\Vert _{\dot{B}{}_{p^{\prime},1}^{-s}}=1}\int_{0}^{\infty}\left\vert \left\langle \tau^{\alpha-1}\nabla\cdot E_{\alpha,\alpha}(-\tau^{\alpha}(-\Delta)^{\theta/2})f(\tau),h\right\rangle \right\vert d\tau\nonumber \\
		& \,\,\,\leq C\mathop{\sup}\limits _{\left\Vert h\right\Vert _{\dot{B}{}_{p^{\prime},1}^{-s}}=1}\int_{0}^{\infty}\left\vert \left\langle f(\tau),\tau^{\alpha-1}\nabla\cdot E_{\alpha,\alpha}(-\tau^{\alpha}(-\Delta)^{\theta/2})h\right\rangle \right\vert d\tau\nonumber \\
		& \,\,\,\leq C\mathop{\sup}\limits _{\left\Vert h\right\Vert _{\dot{B}{}_{p^{\prime},1}^{-s}}=1}\int_{0}^{\infty}\left\Vert f(\tau)\right\Vert _{\dot{B}{}_{p,\infty}^{s_{0}}}\left\Vert \tau^{\alpha-1}\nabla\cdot E_{\alpha,\alpha}(-\tau^{\alpha}(-\Delta)^{\theta/2})h\right\Vert _{\dot{B}{}_{p^{\prime},1}^{-s_{0}}}d\tau\nonumber \\
		& \,\,\,\leq C\mathop{\sup}\limits _{\tau>0}\left\Vert f(\tau)\right\Vert _{\dot{B}{}_{p,\infty}^{s_{0}}}\mathop{\sup}\limits _{\left\Vert h\right\Vert _{\dot{B}{}_{p^{\prime},1}^{-s}}=1}\int_{0}^{\infty}\left\Vert \tau^{\alpha-1}\nabla\cdot E_{\alpha,\alpha}(-\tau^{\alpha}(-\Delta)^{\theta/2})h\right\Vert _{\dot{B}{}_{p^{\prime},1}^{-s_{0}}}d\tau\nonumber \\
		& \,\,\,\le C\mathop{\sup}\limits _{\tau>0}\left\Vert f(\tau)\right\Vert _{\dot{B}{}_{p,\infty}^{s_{0}}}\mathop{\sup}\limits _{\left\Vert h\right\Vert _{\dot{B}{}_{p^{\prime},1}^{-s}}=1}\left\Vert h\right\Vert _{\dot{B}{}_{p^{\prime},1}^{-s}}\nonumber \\
		& \,\,\,\le C\mathop{\sup}\limits _{\tau>0}\left\Vert f(\tau)\right\Vert _{\dot{B}{}_{p,\infty}^{s_{0}}}.\label{aux-est-1000}
	\end{align}
\end{proof}

Now we consider the operator $\mathcal{T}$ defined by $\mathcal{T}(\eta):=\int_{0}^{\infty}s^{\alpha-1}E_{\alpha,\alpha}(-s^\alpha((-\Delta)^{{\theta}/{2}}-\gamma))\eta(s)ds.$
We have the following lemma.

\begin{lemma} \label{Lem: Operador integral T} Let $1<p\le q$,
	$\alpha>0$ and $\theta>0$, then 
	\[
	\left\Vert \mathcal{T}\left(\eta\right)\right\Vert _{\dot{B}{}_{q,\infty}^{2-\theta-\theta_{1}+\frac{n}{q}}}\leq C\mathop{\sup}\limits _{t>0}\left\Vert \eta(t)\right\Vert _{\dot{B}{}_{p,\infty}^{2-2\theta-\theta_{1}+\frac{n}{p}}},
	\]
	for all $\eta\in L^{\infty}\left(\left(0,\infty\right);\dot{B}{}_{p,\infty}^{2-2\theta-\theta_{1}+\frac{n}{p}}\right)$.
	
\end{lemma} \vspace{0.3cm}

\begin{proof} Using duality and Lemma \ref{Lem: Estimativa de integral en Besov Para operador Bilineal}
	with $\zeta=0$ we have
	
	\begin{align}
		& \left\Vert \mathcal{T}\left(\eta\right)\right\Vert _{\dot{B}{}_{q,\infty}^{2-\theta-\theta_{1}+\frac{n}{q}}}=\mathop{\sup}\limits _{\left\Vert h\right\Vert _{\dot{B}{}_{q^{\prime},1}^{-(2-\theta-\theta_{1}+\frac{n}{q})}}=1}\left\vert \left\langle \mathcal{T}\left(\eta\right),h\right\rangle \right\vert \nonumber \\
		& \,\,\,=\mathop{\sup}\limits _{\left\Vert h\right\Vert _{\dot{B}{}_{q^{\prime},1}^{-(2-\theta-\theta_{1}+\frac{n}{q})}}=1}\int_{0}^{\infty}\left\vert \left\langle \tau^{\alpha-1}E_{\alpha,\alpha}(-\tau^{\alpha}(-\Delta)^{\theta/2})\eta(\tau),h\right\rangle \right\vert d\tau\nonumber \\
		& \,\,\,=\mathop{\sup}\limits _{\left\Vert h\right\Vert _{\dot{B}{}_{q^{\prime},1}^{-(2-\theta-\theta_{1}+\frac{n}{q})}}=1}\int_{0}^{\infty}\left\vert \left\langle \eta(\tau),\tau^{\alpha-1}E_{\alpha,\alpha}(-\tau^{\alpha}(-\Delta)^{\theta/2})h\right\rangle \right\vert d\tau\nonumber \\
		& \,\,\,\leq C\mathop{\sup}\limits _{\left\Vert h\right\Vert _{\dot{B}{}_{q^{\prime},1}^{-(2-\theta-\theta_{1}+\frac{n}{q})}}=1}\int_{0}^{\infty}\left\Vert \eta(\tau)\right\Vert _{\dot{B}{}_{p,\infty}^{2-2\theta-\theta_{1}+\frac{n}{p}}}\left\Vert \tau^{\alpha-1}E_{\alpha,\alpha}(-\tau^{\alpha}(-\Delta)^{\theta/2})h\right\Vert _{\dot{B}{}_{p^{\prime},1}^{-\left(2-2\theta-\theta_{1}+\frac{n}{p}\right)}}d\tau\nonumber \\
		& \,\,\,\leq C\mathop{\sup}\limits _{\tau>0}\left\Vert \eta(\tau)\right\Vert _{\dot{B}{}_{p,\infty}^{2-2\theta-\theta_{1}+\frac{n}{p}}}\mathop{\sup}\limits _{\left\Vert h\right\Vert _{\dot{B}{}_{q^{\prime},1}^{-(2-\theta-\theta_{1}+\frac{n}{q})}}=1}\int_{0}^{\infty}\left\Vert \tau^{\alpha-1}E_{\alpha}(-\tau^{\alpha}(-\Delta)^{\theta/2})h\right\Vert _{\dot{B}{}_{p^{\prime},1}^{-\left(2-2\theta-\theta_{1}+\frac{n}{p}\right)}}d\tau\nonumber \\
		& \,\,\,\leq C\mathop{\sup}\limits _{\tau>0}\left\Vert \eta(\tau)\right\Vert _{\dot{B}{}_{p,\infty}^{2-2\theta-\theta_{1}+\frac{n}{p}}}\mathop{\sup}\limits _{\left\Vert h\right\Vert _{\dot{B}{}_{q^{\prime},1}^{-(2-\theta-\theta_{1}+\frac{n}{q})}}=1}\left\Vert h\right\Vert _{\dot{B}{}_{p^{\prime},1}^{-(2-\theta-\theta_{1}+\frac{n}{p})}}\nonumber \\
		& \,\,\,\le C\mathop{\sup}\limits _{\tau>0}\left\Vert \eta(\tau)\right\Vert _{\dot{B}{}_{p,\infty}^{2-2\theta-\theta_{1}+\frac{n}{p}}}\mathop{\sup}\limits _{\left\Vert h\right\Vert _{\dot{B}{}_{q^{\prime},1}^{-(2-\theta-\theta_{1}+\frac{n}{q})}}=1}\left\Vert h\right\Vert _{\dot{B}{}_{q^{\prime},1}^{-(2-\theta-\theta_{1}+\frac{n}{q})}}\nonumber \\
		& \,\,\,\le C\mathop{\sup}\limits _{\tau>0}\left\Vert \eta(\tau)\right\Vert _{\dot{B}{}_{p,\infty}^{2-2\theta-\theta_{1}+\frac{n}{p}}}.\label{eq: Aux operador integral lineal}
	\end{align}
\end{proof}

\begin{remark}
	In the previous proof we have used that $\dot{B}{}_{q^{\prime},1}^{-(2-\theta-\theta_{1}+\frac{n}{q})}\hookrightarrow\dot{B}{}_{p^{\prime},1}^{-(2-\theta-\theta_{1}+\frac{n}{p})}$ which is a direct consequence of Lemma \ref{lem:Bernstein inequality}.
\end{remark}

\subsection{ Product Estimate. Proof of Theorem \ref{lem: produto}.}

%

\begin{proof} For this proof denote $s_{1}=2-2\theta-\theta_{1}+\frac{n}{p}$,
	$s_{0}=3-3\theta-\theta_{1}+\frac{n}{p}$ and $s_{2}=2-\theta-\theta_{1}+\frac{n}{q}.$ 
	
	From the decomposition (\ref{eq: descomposion de Bony}), we obtain
	\begin{align}
		\Delta_{j}\left(fG(g)\right) & =\sum\limits _{\left\vert k-j\right\vert \leq4}\Delta_{j}\left(S_{k-2}f\Delta_{k}G(g)\right)+\sum\limits _{\left\vert k-j\right\vert \leq4}\Delta_{j}\left(S_{k-2}G(g)\Delta_{k}f\right)+\sum\limits _{k\geq j-2}\Delta_{j}\left(\Delta_{k}f\tilde{\Delta}_{k}G(g)\right)\nonumber \\
		& =I_{1}^{j}+I_{2}^{j}+I_{3}^{j}.\label{eq:aux descomposicion paraproducto}
	\end{align}
	In order to estimate $I_{1}^{j}$, let $p^{*}$ such that $\frac{1}{p}=\frac{1}{p^{*}}+\frac{1}{q}$,
	then
	
	\begin{align*}
		& \left\Vert I_{1}^{j}\right\Vert _{L^{p}}\leq C\sum\limits _{\left\vert k-j\right\vert \leq4}\left\Vert S_{k-2}f\right\Vert _{L^{p^{*}}}\left\Vert \Delta_{k}G(g)\right\Vert _{L^{q}}\leq C\sum\limits _{\left\vert k-j\right\vert \leq4}\left(\sum\limits _{m\leq k-2}\left\Vert \Delta_{m}f\right\Vert _{L^{p^{*}}}\right)\left\Vert G(\Delta_{k}g)\right\Vert _{L^{q}}\\
		& \,\,\,\leq C\sum\limits _{\left\vert k-j\right\vert \leq4}\left(\sum\limits _{m\leq k-2}2^{m(\frac{n}{p}-\frac{n}{p^{*}})}\left\Vert \Delta_{m}f\right\Vert _{L^{p}}\right)2^{k\left(1-\theta_{1}\right)}\left\Vert \Delta_{k}g\right\Vert _{L^{q}}\\
		& \,\,\,\leq C\left\Vert f\right\Vert _{\dot{B}{}_{p,\infty}^{s_{1}+\rho_{1}}}\left\Vert g\right\Vert _{\dot{B}{}_{q,\infty}^{s_{2}+\rho_{2}}}\sum\limits _{\left\vert k-j\right\vert \leq4}\left(\sum\limits _{m\leq k-2}2^{m(-2+2\theta+\theta_{1}-\frac{n}{p^{*}}-\rho_{1})}\right)2^{k\left(1-\theta_{1}-2+\theta+\theta_{1}-\frac{n}{q}-\rho_{2}\right)}\\
		& \,\,\,\leq C\left\Vert f\right\Vert _{\dot{B}{}_{p,\infty}^{s_{1}+\rho_{1}}}\left\Vert g\right\Vert _{\dot{B}{}_{q,\infty}^{s_{2}+\rho_{2}}}\sum\limits _{\left\vert k-j\right\vert \leq4}\left(\sum\limits _{m\leq k-2}2^{m(-2+2\theta+\theta_{1}+\frac{n}{q}-\frac{n}{p}-\rho_{1})}\right)2^{k\left(-1+\theta-\frac{n}{q}-\rho_{2}\right)}.
	\end{align*}
	Note that $-2+2\theta+\theta_{1}+\frac{n}{q}-\frac{n}{p}-\rho_{1}>0$
	for some $\rho_{1}$ small enough, in fact, since $\frac{n}{q}\geq\frac{n}{p^{\prime}}$
	we only need to verify that $-2+2\theta+\theta_{1}+\frac{n}{p^{\prime}}-\frac{n}{p}>0,$
	which reduces to the condition $\theta>\ensuremath{1-\frac{n}{2}-\frac{\theta_{1}}{2}+\frac{n}{p}}$.
	Therefore,
	
	\begin{align*}
		& \left\Vert I_{1}^{j}\right\Vert _{L^{p}}\le C\left\Vert f\right\Vert _{\dot{B}{}_{p,\infty}^{s_{1}+\rho_{1}}}\left\Vert g\right\Vert _{\dot{B}{}_{q,\infty}^{s_{2}+\rho_{2}}}\sum\limits _{\left\vert k-j\right\vert \leq4}2^{k(-2+2\theta+\theta_{1}+\frac{n}{q}-\frac{n}{p}-\rho_{1})}2^{k\left(-1+\theta-\frac{n}{q}-\rho_{2}\right)}\\
		& \,\,\,\leq C\left\Vert f\right\Vert _{\dot{B}{}_{p,\infty}^{s_{1}+\rho_{1}}}\left\Vert g\right\Vert _{\dot{B}{}_{q,\infty}^{s_{2}+\rho_{2}}}2^{j\left(-2+2\theta+\theta_{1}+\frac{n}{q}-\frac{n}{p}-\rho_{1}-1+\theta-\frac{n}{q}-\rho_{2}\right)}\\
		& \,\,\,\le C\left\Vert f\right\Vert _{\dot{B}{}_{p,\infty}^{s_{1}+\rho_{1}}}\left\Vert g\right\Vert _{\dot{B}{}_{q,\infty}^{s_{2}+\rho_{2}}}2^{j\left(-3+3\theta+\theta_{1}-\frac{n}{p}-\rho_{1}-\rho_{2}\right)},
	\end{align*}
	and thus, 
	\begin{equation}
		\left\Vert I_{1}^{j}\right\Vert _{L^{p}}\leq C\left\Vert f\right\Vert _{\dot{B}{}_{p,\infty}^{s_{1}+\rho_{1}}}\left\Vert g\right\Vert _{\dot{B}{}_{q,\infty}^{s_{2}+\rho_{2}}}2^{-j\left(s_{0}+\rho_{1}+\rho_{2}\right)}.\label{eq:aux produnto I1}
	\end{equation}
	In order to estimate $I_{2}^{j}$, we proceed similarly to obtain
	\begin{align*}
		& \left\Vert I_{2}^{j}\right\Vert _{L^{p}}\leq C\sum\limits _{\left\vert k-j\right\vert \leq4}\left\Vert S_{k-2}G(g)\right\Vert _{L^{\infty}}\left\Vert \Delta_{k}f\right\Vert _{L^{p}}\leq C\sum\limits _{\left\vert k-j\right\vert \leq4}\left(\sum\limits _{m\leq k-2}\left\Vert \Delta_{m}G(g)\right\Vert _{L^{\infty}}\right)\left\Vert \Delta_{k}f\right\Vert _{L^{p}}\\
		& \,\,\,\leq C\sum\limits _{\left\vert k-j\right\vert \leq4}\left(\sum\limits _{m\leq k-2}2^{m\left(\frac{n}{q}\right)}\left\Vert \Delta_{m}G(g)\right\Vert _{L^{q}}\right)\left\Vert \Delta_{k}f\right\Vert _{L^{p}}\\
		& \,\,\,\leq C\sum\limits _{\left\vert k-j\right\vert \leq4}\left(\sum\limits _{m\leq k-2}2^{m\left(\frac{n}{q}+1-\theta_{1}\right)}\left\Vert \Delta_{m}g\right\Vert _{L^{q}}\right)\left\Vert \Delta_{k}f\right\Vert _{L^{p}}\\
		& \,\,\,\le C\left\Vert g\right\Vert _{\dot{B}{}_{q,\infty}^{s_{2}+\rho_{2}}}\sum\limits _{\left\vert k-j\right\vert \leq4}\left(\sum\limits _{m\leq k-2}2^{m\left(\frac{n}{q}+1-\theta_{1}-2+\theta+\theta_{1}-\frac{n}{q}-\rho_{2}\right)}\right)\left\Vert \Delta_{k}f\right\Vert _{L^{p}}\\
		& \,\,\,\le C\left\Vert f\right\Vert _{\dot{B}{}_{p,\infty}^{s_{1}+\rho_{1}}}\left\Vert g\right\Vert _{\dot{B}{}_{q,\infty}^{s_{2}+\rho_{2}}}\sum\limits _{\left\vert k-j\right\vert \leq4}\left(\sum\limits _{m\leq k-2}2^{m\left(-1+\theta-\rho_{2}\right)}\right)2^{k\left(-2+2\theta+\theta_{1}-\frac{n}{p}\right)}\\
		& \,\,\,\le C\left\Vert f\right\Vert _{\dot{B}{}_{p,\infty}^{s_{1}+\rho_{1}}}\left\Vert g\right\Vert _{\dot{B}{}_{q,\infty}^{s_{2}+\rho_{2}}}\sum\limits _{\left\vert k-j\right\vert \leq4}2^{k\left(-1+\theta-\rho_{2}\right)}2^{k\left(-2+2\theta+\theta_{1}-\frac{n}{p}\right)}\\
		& \,\,\,\le C\left\Vert f\right\Vert _{\dot{B}{}_{p,\infty}^{s_{1}+\rho_{1}}}\left\Vert g\right\Vert _{\dot{B}{}_{q,\infty}^{s_{2}+\rho_{2}}}2^{j\left(-3+3\theta+\theta_{1}-\frac{n}{p}-\rho_{1}-\rho_{2}\right)}.
	\end{align*}
	Here we use that $\theta>1$ which implies that $-1+\theta>0$ and
	so $-1+\theta-\rho_{2}>0$ for $\rho_{2}$ small enough. The previous
	inequality reduces to 
	\begin{equation}
		\left\Vert I_{2}^{j}\right\Vert _{L^{p}}\leq C\left\Vert f\right\Vert _{\dot{B}{}_{p,\infty}^{s_{1}+\rho_{1}}}\left\Vert g\right\Vert _{\dot{B}{}_{q,\infty}^{s_{2}+\rho_{2}}}2^{-j\left(s_{0}+\rho_{1}+\rho_{2}\right)}.\label{eq:aux produnto I2}
	\end{equation}
	
	Now we turn to $I_{3}^{j}$. Note that in the given conditions we
	have $-3+3\theta+\theta_{1}-n<0$ and so $-3+3\theta+\theta_{1}-n-\rho_{1}-\rho_{2}<0$;
	then we have the estimate 
	\begin{align*}
		\left\Vert I_{3}^{j}\right\Vert _{L^{1}} & \leq C\sum\limits _{k\geq j-2}\left\Vert \Delta_{k}f\tilde{\Delta}_{k}G(g)\right\Vert _{L^{1}}\leq C\sum\limits _{k\geq j-2}\left\Vert \Delta_{k}f\right\Vert _{L^{p}}\left\Vert \tilde{\Delta}_{k}G(g)\right\Vert _{L^{p^{\prime}}}\\
		& \leq C\left\Vert f\right\Vert _{\dot{B}{}_{p,\infty}^{s_{1}+\rho_{1}}}\sum\limits _{k\geq j-2}2^{k(-2+2\theta+\theta_{1}-\frac{n}{p}-\rho_{1})}\left\Vert \tilde{\Delta}_{k}G(g)\right\Vert _{L^{p^{\prime}}}\\
		& \le C\left\Vert f\right\Vert _{\dot{B}{}_{p,\infty}^{s_{1}+\rho_{1}}}\sum\limits _{k\geq j-2}2^{k(-2+2\theta+\theta_{1}-\frac{n}{p}-\rho_{1})}2^{k(1-\theta_{1})}\left\Vert \tilde{\Delta}_{k}g\right\Vert _{L^{p^{\prime}}}\\
		& \le C\left\Vert f\right\Vert _{\dot{B}{}_{p,\infty}^{s_{1}+\rho_{1}}}\sum\limits _{k\geq j-2}2^{k(-2+2\theta+\theta_{1}-\frac{n}{p}-\rho_{1})}2^{k(1-\theta_{1})}2^{k\left(\frac{n}{q}-\frac{n}{p^{\prime}}\right)}\left\Vert \tilde{\Delta}_{k}g\right\Vert _{L^{q}}\\
		& \le C\left\Vert f\right\Vert _{\dot{B}{}_{p,\infty}^{s_{1}+\rho_{1}}}\left\Vert g\right\Vert _{\dot{B}{}_{q,\infty}^{s_{2}+\rho_{2}}}\sum\limits _{k\geq j-2}2^{k(-3+3\theta+\theta_{1}-n-\rho_{1}-\rho_{2})}\\
		& \le C\left\Vert f\right\Vert _{\dot{B}{}_{p,\infty}^{s_{1}+\rho_{1}}}\left\Vert g\right\Vert _{\dot{B}{}_{q,\infty}^{s_{2}+\rho_{2}}}2^{j(-3+3\theta+\theta_{1}-n-\rho_{1}-\rho_{2})}.
	\end{align*}
	So
	
	\begin{align}
		\left\Vert I_{3}^{j}\right\Vert _{L^{p}} & \le C2^{j\left(\frac{n}{1}-\frac{n}{p}\right)}\left\Vert I_{3}^{j}\right\Vert _{L^{1}}\nonumber \\
		& \le C2^{j\left(n-\frac{n}{p}\right)}\left\Vert f\right\Vert _{\dot{B}{}_{p,\infty}^{s_{1}+\rho_{1}}}\left\Vert g\right\Vert _{\dot{B}{}_{q,\infty}^{s_{2}+\rho_{2}}}2^{j(-3+3\theta+\theta_{1}-n-\rho_{1}-\rho_{2})}\nonumber \\
		& \leq C\left\Vert f\right\Vert _{\dot{B}{}_{p,\infty}^{s_{1}+\rho_{1}}}\left\Vert g\right\Vert _{\dot{B}{}_{q,\infty}^{s_{2}+\rho_{2}}}2^{-j\left(s_{0}+\rho_{1}+\rho_{2}\right)}.\label{eq:aux produnto I3}
	\end{align}
	
	Computing the norm $\left\Vert \cdot\right\Vert _{L^{p}}$ in (\ref{eq:aux descomposicion paraproducto})
	and considering the estimates $(\ref{eq:aux produnto I1})$, $(\ref{eq:aux produnto I2})$
	and $(\ref{eq:aux produnto I3})$, we get the result.
\end{proof}

\begin{remark}\label{rem:aux1}
	
	In the proof of Lemma $\ref{eq:EstimativaProducto},$ the condition
	$6n/(\theta_{1}+5n)<p$ is not directly used. This condition is imposed
	in order to guarantee that the interval $\left(\ensuremath{1-\frac{n}{2}-\frac{\theta_{1}}{2}+\frac{n}{p}},1+\frac{n-\theta_{1}}{3}\right)$
	for $\theta$ is nonempty. Also, the condition $\theta_{1}\in[0,n)$
	guarantees that the interval $\left(1,1+\frac{n-\theta_{1}}{3}\right)$
	for $\theta$ is nonempty.
	
\end{remark}

\subsection{Bilinear estimate. Proof of Theorem \ref{Bilinearfinal}.}

\begin{proof} Let $0<T\leq\infty$ and $t\in\left(0,T\right).$
	The bilinear term $B_{\theta}(\eta,v)$ can be written as
	
	\[
	B_{\theta}(\eta,v)(t)=-\int_{0}^{t}(t-\tau)^{\alpha-1}\nabla\cdot E_{\alpha,\alpha}(-(t-\tau)^{\alpha}(-\Delta)^{\theta/2})(\eta G(v))d\tau=\mathcal{B}(f_{t}),
	\]
	where $f_{t}(x,\tau)$ is defined by 
	\begin{align*}
		f_{t}(\cdot,\tau) & =(\eta G(v))\left(\cdot,t-\tau\right),\,\,\mbox{}\,\mbox{a.e.}\,\tau\in\left(0,t\right),\\
		f_{t}(\cdot,\tau) & =0,\,\,\mbox{\,\ a.e.\,}\tau\in\left(t,\infty\right).
	\end{align*}
	From Lemma \ref{Lem: Operador integral B} (with $s=2-2\theta-\theta_{1}+\frac{n}{p}$
	and $s_{0}=3-3\theta-\theta_{1}+\frac{n}{p}$) we get 
	\[
	\left\Vert \mathcal{B}\left(f_{t}\right)\right\Vert _{\dot{B}{}_{p,\infty}^{2-2\theta-\theta_{1}+\frac{n}{p}}}\leq C\mathop{\sup}\limits _{\tau>0}\left\Vert f_{t}(\tau)\right\Vert _{\dot{B}{}_{p,\infty}^{3-3\theta-\theta_{1}+\frac{n}{p}}}.
	\]
	Using Lemma \ref{lem: produto}$,$ with $\rho_{1}=\rho_{2}=0,$ we
	can estimate
	\begin{align*}
		\mathop{\sup}\limits _{0<\tau<T}\left\Vert f_{t}(\tau)\right\Vert _{\dot{B}{}_{p,\infty}^{3-3\theta-\theta_{1}+\frac{n}{p}}} & =\mathop{\sup}\limits _{0<\tau<t<T}\left\Vert (\eta G(v))\left(\cdot,t-\tau\right)\right\Vert _{\dot{B}{}_{p,\infty}^{3-3\theta-\theta_{1}+\frac{n}{p}}}\\
		& \leq C\mathop{\sup}\limits _{0<\tau<t<T}\left\Vert \eta\left(\cdot,t-\tau\right)\right\Vert _{\dot{B}{}_{p,\infty}^{2-2\theta-\theta_{1}+\frac{n}{p}}}\left\Vert v\left(\cdot,t-\tau\right)\right\Vert _{\dot{B}{}_{q,\infty}^{2-\theta-\theta_{1}+\frac{n}{q}}}\\
		& \leq C\mathop{\sup}\limits _{0<\tau<t<T}\left\Vert \eta\left(\cdot,t-\tau\right)\right\Vert _{\dot{B}{}_{p,\infty}^{2-2\theta-\theta_{1}+\frac{n}{p}}}\mathop{\sup}\limits _{0<s<t<T}\left\Vert v\left(\cdot,t-\tau\right)\right\Vert _{\dot{B}{}_{q,\infty}^{2-\theta-\theta_{1}+\frac{n}{q}}}\\
		& \le C\mathop{\sup}\limits _{0<\tau<T}\left\Vert \eta\left(\cdot,\tau\right)\right\Vert _{\dot{B}{}_{p,\infty}^{2-2\theta-\theta_{1}+\frac{n}{p}}}\mathop{\sup}\limits _{0<\tau<T}\left\Vert v\left(\cdot,\tau\right)\right\Vert _{\dot{B}{}_{q,\infty}^{2-\theta-\theta_{1}+\frac{n}{q}}}.
	\end{align*}
	Thus, we can conclude that 
	\[
	\mathop{\sup}\limits _{0<t<T}\left\Vert B_{\theta}(\eta,v)(t)\right\Vert _{\dot{B}{}_{p,\infty}^{2-2\theta-\theta_{1}+\frac{n}{p}}}\leq K\mathop{\sup}\limits _{0<t<T}\left\Vert \eta\left(t\right)\right\Vert _{\dot{B}{}_{p,\infty}^{2-2\theta-\theta_{1}+\frac{n}{p}}}\mathop{\sup}\limits _{0<t<T}\left\Vert v\left(t\right)\right\Vert _{\dot{B}{}_{q,\infty}^{2-\theta-\theta_{1}+\frac{n}{q}}}.
	\]
\end{proof}

\begin{lemma}{\label{linearfinal}}
	
	Let $1<p\le q$, $\alpha>0$ and $\theta>0$. Then, there exists a
	constant $K>0$ (independent of $T$) such that
	
	\[
	\left\Vert T_{\theta}\left(\eta\right)\right\Vert _{L^{\infty}\left((0,T);\dot{B}{}_{q,\infty}^{2-\theta-\theta_{1}+\frac{n}{q}}\right)}\leq C\left\Vert \eta(t)\right\Vert _{L^{\infty}\left((0,T);\dot{B}{}_{p,\infty}^{2-2\theta-\theta_{1}+\frac{n}{p}}\right)},
	\]
	for all $\eta\in L^{\infty}\left((0,T);\dot{B}{}_{p,\infty}^{s}\right)$.
	
\end{lemma}

\begin{proof} Let $0<T\leq\infty$ and $t\in\left(0,T\right).$
	Note that the operator $T(\eta)$ can be written as
	
	\[
	T_{\theta}(\eta)=\int_{0}^{t}(t-\tau)^{\alpha-1}E_{\alpha,\alpha}(-(t-\tau)^{\alpha}((-\Delta)^{\theta/2}-\gamma))\eta(\tau)d\tau=-\mathcal{T}(f_{t}),
	\]
	where $f_{t}(x,\tau)$ is defined by 
	\begin{align*}
		f_{t}(\cdot,\tau) & =\eta\left(\cdot,t-\tau\right),\,\,\mbox{}\,\mbox{a.e.}\,\tau\in\left(0,t\right),\\
		f_{t}(\cdot,\tau) & =0,\,\,\mbox{\,\ a.e.\,}\tau\in\left(t,\infty\right).
	\end{align*}
	From Lemma \ref{Lem: Operador integral T} we obtain
	
	\[
	\left\Vert \mathcal{T}\left(f_{t}\right)\right\Vert _{\dot{B}{}_{q,\infty}^{2-\theta-\theta_{1}+\frac{n}{q}}}\leq C\mathop{\sup}\limits _{\tau>0}\left\Vert f_{t}(\tau)\right\Vert _{\dot{B}{}_{p,\infty}^{2-2\theta-\theta_{1}+\frac{n}{p}}},
	\]
	moreover 
	\begin{align*}
		\mathop{\sup}\limits _{0<\tau<T}\left\Vert f_{t}(\tau)\right\Vert _{\dot{B}{}_{p,\infty}^{2-2\theta-\theta_{1}+\frac{n}{p}}} & =\mathop{\sup}\limits _{0<\tau<t<T}\left\Vert \eta\left(t-\tau\right)\right\Vert _{\dot{B}{}_{p,\infty}^{2-2\theta-\theta_{1}+\frac{n}{p}}}\\
		& \le\mathop{\sup}\limits _{0<\tau<T}\left\Vert \eta\left(\cdot,\tau\right)\right\Vert _{\dot{B}{}_{p,\infty}^{2-2\theta-\theta_{1}+\frac{n}{p}}}.
	\end{align*}
	Thus, we arrive at 
	\[
	\mathop{\sup}\limits _{0<t<T}\left\Vert T_{\theta}(\eta)\right\Vert _{\dot{B}{}_{q,\infty}^{2-\theta-\theta_{1}+\frac{n}{q}}}\leq K\mathop{\sup}\limits _{0<\tau<T}\left\Vert \eta\left(\cdot,\tau\right)\right\Vert _{\dot{B}{}_{p,\infty}^{2-2\theta-\theta_{1}+\frac{n}{p}}}.
	\]
\end{proof}

\section{Existence and uniqueness of global solutions}\label{proofs}
The aim of this section is to prove the existence and uniqueness of global mild solution of system (\ref{eq: Keller-Segel intro}), which will be carried out through an iterative approach.  

\subsection{Proof of Theorem \ref{Teo: Main}.}

To simplify the notation, we denote 
\[
X=L^{\infty}\left((0,T);\dot{B}{}_{p,\infty}^{2-2\theta-\theta_{1}+\frac{n}{p}}\right)\ \mbox{and}\ Y=L^{\infty}\left((0,T);\dot{B}{}_{q,\infty}^{2-\theta-\theta_{1}+\frac{n}{q}}\right).
\]
In order to prove Theorem \ref{Teo: Main}, we consider the following iterative
system: 
\[
\eta^{1}:=E_{\alpha}(-t^\alpha(-\Delta)^{{\theta}/{2}})\eta_{0},\qquad v^{1}:=E_{\alpha}(-t^\alpha((-\Delta)^{{\theta}/{2}}-\gamma))v_{0},
\]
and for $n\geq1$ 
\[
\begin{array}{c}
	\eta^{n+1}:=\eta^{1}+B_\theta(\eta^{n},v^{n}),\\
	v^{n+1}:=v^{1}+T_\theta(\eta^{n+1}).
\end{array}
\]
From Lemma \ref{Lem: Frac_Heat-Kernel-estimate en Besov}, it follows that 
\[
\|\eta^{1}\|_{X}=\|E_{\alpha}(-t^\alpha(-\Delta)^{{\theta}/{2}})\eta_{0}\|_{X}\le C_{1}\|\eta_{0}\|_{X},
\]
and 
\[
\|v^{1}\|_{Y}=\|E_{\alpha}(-t^\alpha((-\Delta)^{{\theta}/{2}}-\gamma))v_{0}\|_{Y}\le C_{2}\|v_{0}\|_{Y}.
\]
Additionally, using Lemmas \ref{Bilinearfinal} and \ref{linearfinal}  we have
\begin{align*}
	\|\eta^{n+1}\|_{X} & =\|\eta^{1}+B_\theta(\eta^{n},v^{n})\|_{X}\le\|\eta^{1}\|_{X}+\|B_\theta(\eta^{n},v^{n})\|_{X}\\
	& \le C_{1}\|\eta_{0}\|_{X}+K\|\eta^{n}\|_{X}\|v^{n}\|_{Y},
\end{align*}
\begin{align*}
	\|v^{n+1}\|_{Y} & =\|v^{1}+T_\theta(\eta^{n+1})\|_{Y}\le\|v^{1}\|_{Y}+\|T_\theta(\eta^{n+1})\|_{Y}\\
	& \le C_{2}\|v_{0}\|_{Y}+C\|\eta^{n+1}\|_{X}.
\end{align*}
Let $0<\varepsilon<\frac{1}{2K}$ and $\eta_{0}$,$v_{0}$ such that
$C_{1}\|\eta_{0}\|_{X}\le\frac{\varepsilon}{4C}<\frac{\varepsilon}{2C}$
and $C_{2}\|v_{0}\|_{Y}\le\frac{\varepsilon}{2}<\varepsilon,$ then
\begin{align*}
	\|\eta^{2}\|_{X} & <\frac{\varepsilon}{4C}+K\frac{\varepsilon}{2C}\varepsilon<\frac{\varepsilon}{4C}+\frac{\varepsilon}{4C}=\frac{\varepsilon}{2C},
\end{align*}
\begin{align*}
	\|v^{2}\|_{Y} & <\frac{\varepsilon}{2}+C\frac{\varepsilon}{2C}=\varepsilon.
\end{align*}
Proceeding inductively, we prove that 
\[
\|\eta^{n+1}\|_{X}<\frac{\varepsilon}{2C}\qquad\text{and}\qquad\|v^{n+1}\|_{Y}<\varepsilon.
\]
Now we prove that the sequences $\left(\eta^{n}\right)$ and $\left(v^{n}\right)$
are Cauchy in the respective spaces. In fact, we have that

\begin{align*}
	\eta^{n+1}-\eta^{n} & =B_\theta(\eta^{n},v^{n})-B_\theta(\eta^{n-1},v^{n-1})\\
	& =B_\theta(\eta^{n}-\eta^{n-1},v^{n})+B_\theta(\eta^{n-1},v^{n}-v^{n-1}),
\end{align*}
and so 
\begin{align}
	\|\eta^{n+1}-\eta^{n}\|_{X} & \le K\left(\|\eta^{n}-\eta^{n-1}\|_{X}\|v^{n}\|_{Y}+\|\eta^{n-1}\|_{X}\|v^{n}-v^{n-1}\|_{Y}\right)\nonumber \\
	& \le K\left(\varepsilon\|\eta^{n}-\eta^{n-1}\|_{X}+\frac{\varepsilon}{2C}\|v^{n}-v^{n-1}\|_{Y}\right).\label{eq:auxCauchyeta}
\end{align}
On the other hand, since

\[
v^{n+1}-v^{n}=T_\theta(\eta^{n+1}-\eta^{n1}),
\]
we have
\begin{equation}
	\|v^{n+1}-v^{n}\|_{Y}\le C\|\eta^{n+1}-\eta^{n}\|_{X}.\label{eq:auxCauchyv}
\end{equation}
Now, using (\ref{eq:auxCauchyv}) in (\ref{eq:auxCauchyeta}) we arrive
at 
\begin{align}
	\|\eta^{n+1}-\eta^{n}\|_{X} & \le K\left(\|\eta^{n}-\eta^{n-1}\|_{X}\varepsilon+\frac{\varepsilon}{2C}C\|\eta^{n}-\eta^{n-1}\|_{X}\right)\nonumber \\
	& \le\frac{3K\varepsilon}{2}\|\eta^{n}-\eta^{n-1}\|_{X}\le C(\varepsilon)\|\eta^{n}-\eta^{n-1}\|_{X}.\label{eq:auxCauchyeta2}
\end{align}
Under an additional condition on $\varepsilon$ (if required) we ensure
that $C(\epsilon)<1$, and follows from (\ref{eq:auxCauchyeta2})
that $\left(\eta^{n}\right)$ is Cauchy. Finally, from (\ref{eq:auxCauchyv})
we also have that $\left(v^{n}\right)$ is Cauchy. Let $\eta$ and $v$ be such that $\eta^{n}\longrightarrow\eta$ and
$v^{n}\longrightarrow v$. We have that
\begin{align}
	0 & \le\|\eta-E_{\alpha}(-t^{\alpha}(-\Delta)^{\theta/2})\eta_{0}-B(\eta,v)\|_{X}=\|\eta-E_{\alpha}(-t^{\alpha}(-\Delta)^{\theta/2})\eta_{0}-B(\eta,v)-\eta^{n}+\eta^{n}\|_{X}\nonumber \\
	& \le\|\eta-\eta^{n}\|_{X}+\|-E_{\alpha}(-t^{\alpha}(-\Delta)^{\theta/2})\eta_{0}-B(\eta,v)+\eta^{n}\|_{X}\nonumber \\
	& \le\|\eta-\eta^{n}\|_{X}+\|-B(\eta,v)+B(\eta^{n-1},v^{n-1})\|_{X}\nonumber \\
	& \le\|\eta-\eta^{n}\|_{X}+\|B(\eta^{n-1}-\eta,v^{n-1})+B(\eta,v^{n-1}-v)\|_{X}\nonumber \\
	& \le\|\eta-\eta^{n}\|_{X}+K\left(\|\eta^{n-1}-\eta\|_{X}\varepsilon+\frac{\varepsilon}{2C}\|v^{n-1}-v\|_{Y}\right)\nonumber \\
	& \longrightarrow0,\label{eq:etasolucion}
\end{align}
and
\begin{align}
	0 & \le\|v-E_{\alpha}(-t^{\alpha}((-\Delta)^{\theta/2}-\gamma))v_{0}-T(\eta)\|_{Y}=\|v-E_{\alpha}(-t^{\alpha}((-\Delta)^{\theta/2}-\gamma))v_{0}-T(\eta)-v^{n}+v^{n}\|_{Y}\nonumber \\
	& \le\|v-v^{n}\|_{Y}+\|-E_{\alpha}(-t^{\alpha}((-\Delta)^{\theta/2}-\gamma))v_{0}-T(\eta)+v^{n}\|_{Y}\nonumber \\
	& \le\|v-v^{n}\|_{Y}+\|T(\eta^{n}-\eta)\|_{Y}\nonumber \\
	& \le\|v-v^{n}\|_{Y}+C\|\eta^{n}-\eta\|_{X}\nonumber \\
	& \longrightarrow0.\label{eq:vsolucion}
\end{align}
The estimates (\ref{eq:etasolucion}) and (\ref{eq:vsolucion}) prove
that $[\eta,v]$ is a mild solution of (\ref{eq: Keller-Segel intro}). Finally, to prove the uniqueness suppose that $[\eta,v]$ and $[\tilde{\eta},\tilde{v}]$
are two solutions in the same conditions of Theorem \ref{Teo: Main}
with the same initial data, then, following the same ideas of the
proof of (\ref{eq:auxCauchyeta2}) we have
\begin{align*}
	\|\eta-\tilde{\eta}\|_{X} & \le C(\varepsilon)\|\eta-\tilde{\eta}\|_{X},
\end{align*}
with $0<C(\varepsilon)<1$, which implies that $\eta-\tilde{\eta}=0$,
this is $\eta=\tilde{\eta}.$ From this is obvious that $v=\tilde{v}$.

\subsection{Proof of Corollary \ref{Cor: Corolario 1}.}

First, note that 
\begin{align*}
	\mathcal{F}(U_{\theta}(\tau t^{\alpha})\eta_{0}(\sigma\cdot))(\xi) & =e^{-\tau t^{\alpha}\left|\xi\right|^{\theta}}\mathcal{F}(\eta_{0}(\sigma\cdot))=\sigma^{-n}e^{-\tau t^{\alpha}\left|\xi\right|^{\theta}}\mathcal{F}(\eta_{0})(\xi/\sigma)\\
	& =\sigma^{-n}e^{-\tau\left(\sigma^{\frac{\theta}{\alpha}}t\right)^{\alpha}\left|\xi/\sigma\right|^{\theta}}\mathcal{F}(\eta_{0})(\xi/\sigma)\\
	& =\mathcal{F}\left(\left(U_{\theta}\left(\tau\left(\sigma^{\frac{\theta}{\alpha}}t\right)^{\alpha}\right)\eta_{0}\right)(\sigma\cdot)\right)(\xi).
\end{align*}
Thus,
\[
U_{\theta}(\tau t^{\alpha})\eta_{0}(\sigma\cdot)=\left(U_{\theta}\left(\tau\left(\sigma^{\frac{\theta}{\alpha}}t\right)^{\alpha}\right)\eta_{0}\right)(\sigma\cdot),
\]
and 
\begin{align*}
	E_{\alpha}(-t^{\alpha}(-\Delta)^{\theta/2})\eta_{0} & =\int_{0}^{\infty}M_{\alpha}(\tau)U_{\theta}(\tau t^{\alpha})\eta_{0}d\tau=\int_{0}^{\infty}M_{\alpha}(\tau)U_{\theta}(\tau t^{\alpha})\eta_{\sigma0}d\tau\\
	& =\sigma^{2\theta+\theta_{1}-2}\int_{0}^{\infty}M_{\alpha}(\tau)U_{\theta}(\tau t^{\alpha})\eta_{0}(\sigma\cdot)d\tau\\
	& =\sigma^{2\theta+\theta_{1}-2}\int_{0}^{\infty}M_{\alpha}(\tau)\left(U_{\theta}\left(\tau\left(\sigma^{\frac{\theta}{\alpha}}t\right)^{\alpha}\right)\eta_{0}\right)(\sigma\cdot)d\tau,
\end{align*}
this is, 
\[
\eta^{1}(x,t)=\sigma^{2\theta+\theta_{1}-2}\eta^{1}(\sigma x,\sigma^{\frac{\theta}{\alpha}}t).
\]

Using induction, is direct to verify that the members of the sequences
$\left(\eta^{n}\right)$ and $\left(v^{n}\right)$ are invariant by
the scaling \eqref{eq: Scaling}, this is 
\[
\eta^{n}(x,t)=\sigma^{2\theta+\theta_{1}-2}\eta^{n}\left(\sigma x,\sigma^{\frac{\theta}{\alpha}}t\right)\qquad\mbox{and}\qquad v^{n}(x,t)=\sigma^{\theta+\theta_{1}-2}v^{n}\left(\sigma x,\sigma^{\frac{\theta}{\alpha}}t\right).
\]
Finally, since the solution $[\eta,v]$ is the limit in $X\times Y$ of the
sequence $[\eta^{n},v^{n}]$ and the spaces $X,Y$ are invariant for
the scaling, we can conclude that
\[
\eta(x,t)=\sigma^{2\theta+\theta_{1}-2}\eta\left(\sigma x,\sigma^{\frac{\theta}{\alpha}}t\right)\qquad\mbox{and}\qquad v(x,t)=\sigma^{\theta+\theta_{1}-2}v\left(\sigma x,\sigma^{\frac{\theta}{\alpha}}t\right),
\]
this is, $[\eta,v]$ is a self-similar solution.

\subsection{Uniqueness. Proof of Theorem \ref{Cor: Unicidadsinnorma}.}

With the bilinear estimate \eqref{eq:EcuacionBilinear} in hands and
the correct use of the product estimate \eqref{eq:EstimativaProducto},
the uniqueness follows by adapting an argument due to Meyer \cite{Meyer}.
For this proof denote $s_{1}=2-2\theta-\theta_{1}+\frac{n}{p}$, $s_{0}=3-3\theta-\theta_{1}+\frac{n}{p}$
and $s_{2}=2-\theta-\theta_{1}+\frac{n}{q},$ and let $[\eta^{1},v^{1}]$
and $[\eta^{2},v^{2}]$ be two mild solutions in $C\left([0,T);\dot{B}{}_{p,\infty}^{s_{1}}\right)\times C\left([0,T);\dot{B}{}_{q,\infty}^{s_{2}}\right)$
with the same initial data $[\eta_{0},v_{0}]$ in $\tilde{\dot{B}}{}_{p,\infty}^{s_{1}}\times\tilde{\dot{B}}{}_{q,\infty}^{s_{2}}$.
First we prove that there exists $0<T_{1}<T$ such that $[\eta^{1}(t),v^{1}(t)]=[\eta^{2}(t),v^{2}(t)]$
in $\dot{B}{}_{p,\infty}^{s_{1}}\times\dot{B}{}_{q,\infty}^{s_{2}}$
for all $\ensuremath{t\in\left[0,T_{1}\right)}$. Denoting 
\begin{eqnarray*}
	N & = & \eta^{1}-\eta^{2},V=v^{1}-v^{2},\\
	N_{1} & = & E_{\alpha}(-t^{\alpha}((-\Delta)^{\theta/2}))\eta_{0}-\eta^{1},\\
	N_{2} & = & E_{\alpha}(-t^{\alpha}((-\Delta)^{\theta/2}))\eta_{0}-\eta^{2},\\
	V_{1} & = & E_{\alpha}(-t^{\alpha}((-\Delta)^{\theta/2}-\gamma))v_{0}-v^{1},\\
	V_{2} & = & E_{\alpha}(-t^{\alpha}((-\Delta)^{\theta/2}-\gamma))(t)v_{0}-v^{2},
\end{eqnarray*}
we have that 
\begin{equation}
	V=T_{\theta}(N).\label{eq:Relaci=00003D00003D0000F3n V y N}
\end{equation}
Then, 
\begin{equation}
	\mathop{\sup}\limits _{0<t<T_{1}}\left\Vert V\right\Vert _{\dot{B}{}_{q,\infty}^{s_{2}}}\leq K\mathop{\sup}\limits _{0<t<T_{1}}\left\Vert N\right\Vert _{\dot{B}{}_{p,\infty}^{s_{1}}}.\label{eq:Aux Relacion V y N}
\end{equation}
Moreover, 
\begin{align*}
	& \eta^{1}G(v^{1})-\eta^{2}G(v^{2})=NG(v^{1})+\eta^{2}G(V)\\
	& \ \ \ \ \ =NG\left(E_{\alpha}(-t^{\alpha}((-\Delta)^{\theta/2}-\gamma))v_{0}-V_{1}\right)+\left(E_{\alpha}(-t^{\alpha}((-\Delta)^{\theta/2}))\eta_{0}-N_{2}\right)G(V)\\
	& \ \ \ \ \ =NG\left(E_{\alpha}(-t^{\alpha}((-\Delta)^{\theta/2}-\gamma))v_{0}\right)+\left(E_{\alpha}(-t^{\alpha}((-\Delta)^{\theta/2}))\eta_{0}\right)G(V)-NG(V_{1})-N_{2}G(V).
\end{align*}
Thus, 
\begin{align*}
	\|N\|_{\dot{B}{}_{p,\infty}^{s_{1}}} & =\|B(\eta^{1},v^{1})-B(\eta^{2},v^{2})\|_{\dot{B}{}_{p,\infty}^{s_{1}}}\\
	& =\left\Vert \int_{0}^{t}(t-\tau)^{\alpha-1}\nabla\cdot E_{\alpha,\alpha}(-(t-\tau)^{\alpha}((-\Delta)^{\theta/2}))\left((\eta^{1}G(v^{1}))-(\eta^{2}G(v^{2}))\right)(\tau)d\tau\right\Vert _{\dot{B}{}_{p,\infty}^{s_{1}}}\\
	& \le\left\Vert \int_{0}^{t}(t-\tau)^{\alpha-1}\nabla\cdot E_{\alpha,\alpha}(-(t-\tau)^{\alpha}((-\Delta)^{\theta/2}))\left(NG(V_{1})+N_{2}G(V)\right)(\tau)d\tau\right\Vert _{\dot{B}{}_{p,\infty}^{s_{1}}}\\
	& \,\,\,\,\,+\left\Vert \int_{0}^{t}(t-\tau)^{\alpha-1}\nabla\cdot E_{\alpha,\alpha}(-(t-\tau)^{\alpha}((-\Delta)^{\theta/2}))(NG\left(E_{\alpha}(-\tau^{\alpha}((-\Delta)^{\theta/2}-\gamma))v_{0}\right)\right.\\
	& \,\,\,\,\,+\left.\left(E_{\alpha}(-\tau^{\alpha}((-\Delta)^{\theta/2}))\eta_{0}\right)G(V))d\tau\right\Vert _{\dot{B}{}_{p,\infty}^{s_{1}}}\\
	& \ensuremath{:=J_{1}(t)+J_{2}(t)}.
\end{align*}
For $J_{1}(t)$ we have 
\[
J_{1}(t)\le K\left(\ensuremath{\sup_{0<t<T_{1}}\|N\|_{\dot{B}{}_{p,\infty}^{s_{1}}}}\sup_{0<t<T_{1}}\|V_{1}\|_{\dot{B}{}_{q,\infty}^{s_{2}}}+\ensuremath{\sup_{0<t<T_{1}}\|N_{2}\|_{\dot{B}{}_{p,\infty}^{s_{1}}}}\sup_{0<t<T_{1}}\|V\|_{\dot{B}{}_{q,\infty}^{s_{2}}}\right).
\]
Using $(\ref{eq:Aux Relacion V y N})$ we arrive at

\begin{equation}
	J_{1}(t)\le C\ensuremath{\sup_{0<t<T_{1}}\|N\|_{\dot{B}{}_{p,\infty}^{s_{1}}}}\left(\sup_{0<t<T_{1}}\|V_{1}\|_{\dot{B}{}_{q,\infty}^{s_{2}}}+\ensuremath{\sup_{0<t<T_{1}}\|N_{2}\|_{\dot{B}{}_{p,\infty}^{s_{1}}}}\right).\label{eq:EstimativaJ1}
\end{equation}
On the other hand, for $J_{2}(t)$ it follows that 
\begin{align}
	J_{2}(t) & =\left\Vert \int_{0}^{t}(t-\tau)^{\alpha-1}\nabla\cdot E_{\alpha,\alpha}(-(t-\tau)^{\alpha}((-\Delta)^{\theta/2}))(NG\left(E_{\alpha}(-\tau^{\alpha}((-\Delta)^{\theta/2}-\gamma))v_{0}\right)\right.\nonumber \\
	& \left.\left.\ \ \ \ +\left(E_{\alpha}(-\tau^{\alpha}(-\Delta)^{\theta/2})\eta_{0}\right)G(V)\right)d\tau\right\Vert _{\dot{B}{}_{p,\infty}^{s_{1}}}\nonumber \\
	& \le\int_{0}^{t}\|(t-\tau)^{\alpha-1}\nabla\cdot E_{\alpha,\alpha}(-(t-\tau)^{\alpha}(-\Delta)^{\theta/2})\left(NG\left(E_{\alpha}(-\tau^{\alpha}((-\Delta)^{\theta/2}-\gamma))v_{0}\right)\right.\nonumber \\
	& \ \ \ \ \left.+\left(E_{\alpha}(-\tau^{\alpha}(-\Delta)^{\theta/2})\eta_{0}\right)G(V)\right)\|_{\dot{B}{}_{p,\infty}^{s_{1}}}d\tau\nonumber \\
	& \le C\int_{0}^{t}(t-\tau)^{\alpha-1}(t-\tau)^{-\frac{\alpha}{\theta}\left(s-(s_{0}+\rho)+1\right)}\left\Vert NG\left(E_{\alpha}(-\tau^{\alpha}((-\Delta)^{\theta/2}-\gamma))v_{0}\right)\right.\nonumber \\
	& \ \ \ \ \left.+\left(E_{\alpha}(-\tau^{\alpha}(-\Delta)^{\theta/2})\eta_{0}\right)G(V)\right\Vert _{\dot{B}{}_{p,\infty}^{s_{0}+\rho}}d\tau\nonumber \\
	& \le C\int_{0}^{t}\!(t-\tau)^{\alpha-\!1\!-\frac{\alpha}{\theta}\left(\theta\!-\!\rho\right)}\left\Vert \!NG\left(E_{\alpha}(-\tau^{\alpha}((-\Delta)^{\theta/2}\!-\!\gamma))v_{0}\right)\!+\!\left(E_{\alpha}(-\tau^{\alpha}(-\Delta)^{\theta/2})\eta_{0}\right)G(V)\right\Vert _{\dot{B}{}_{p,\infty}^{s_{0}+\rho}}d\tau\nonumber \\
	& \le C\int_{0}^{t}(t-\tau)^{-1+\frac{\alpha}{\theta}\rho}\left\Vert NG\left(E_{\alpha}(-\tau^{\alpha}((-\Delta)^{\theta/2}-\gamma))v_{0}\right)+\left(E_{\alpha}(-\tau^{\alpha}(-\Delta)^{\theta/2})\eta_{0}\right)G(V)\right\Vert _{\dot{B}{}_{p,\infty}^{s_{0}+\rho}}d\tau\nonumber \\
	& \le C\int_{0}^{t}(t-\tau)^{-1+\frac{\alpha}{\theta}\rho}\tau^{-\frac{\alpha}{\theta}\rho}\!\left\Vert N\right\Vert _{\dot{B}{}_{p,\infty}^{s_{1}}}\left(\tau^{\frac{\alpha}{\theta}\rho}\left\Vert E_{\alpha}(-\tau^{\alpha}((-\Delta)^{\theta/2}-\gamma))v_{0}\right\Vert _{\dot{B}{}_{q,\infty}^{s_{2}+\rho}}\right)d\tau\nonumber \\
	& \ \ \ \ +C\int_{0}^{t}(t-\tau)^{-1+\frac{\alpha}{\theta}\rho}\tau^{\frac{\alpha}{\theta}\rho}\!\left\Vert N\right\Vert _{\dot{B}{}_{p,\infty}^{s_{1}}}\left\Vert E_{\alpha}(-\tau^{\alpha}(-\Delta)^{\theta/2})\eta_{0}\right\Vert _{\dot{B}{}_{p,\infty}^{s_{1}+\rho}}d\tau\nonumber \\
	& \le\left(\sup_{0<t<T_{1}}t^{\frac{\alpha}{\theta}\rho}\left\Vert E_{\alpha}(-t^{\alpha}((-\Delta)^{\theta/2}-\gamma))v_{0}\right\Vert _{\dot{B}{}_{q,\infty}^{s_{2}+\rho}}+\sup_{0<t<T_{1}}t^{\frac{\alpha}{\theta}\rho}\left\Vert E_{\alpha}(-t^{\alpha}(-\Delta)^{\theta/2})\eta_{0}\right\Vert _{\dot{B}{}_{p,\infty}^{s_{1}+\rho}}\right)\nonumber \\
	& \ \ \ \ \times C\sup_{0<t<T_{1}}\left\Vert N\right\Vert _{\dot{B}{}_{p,\infty}^{s_{1}}},\label{eq:EstimativaJ2}
\end{align}
where we have used Lemma $\ref{lem: produto}$ adequately, the relation
$(\ref{eq:Aux Relacion V y N})$, and the fact that for $\rho>0$
small enough the following integral holds 
\[
\int_{0}^{t}(t-\tau)^{-1+\frac{\alpha}{\theta}\rho}\tau^{-\frac{\alpha}{\theta}\rho}d\tau=\int_{0}^{1}(1-r)^{-1+(1-\alpha)+\frac{\alpha}{\theta}\rho}r^{-\frac{\alpha}{\theta}\rho}dr=C.
\]
Thus, from the estimates $\eqref{eq:EstimativaJ1}$ and $\eqref{eq:EstimativaJ2}$
we get 
\begin{align}
	\sup_{0<t<T_{1}}\|N\|_{\dot{B}{}_{p,\infty}^{s_{1}}} & \le CZ(T_{1})\sup_{0<t<T_{1}}\left\Vert N\right\Vert _{\dot{B}{}_{p,\infty}^{s_{1}}},\label{eq: Estimativa final eta}
\end{align}
with 
\begin{align*}
	Z(T_{1})= & \sup_{0<t<T_{1}}\|V_{1}\|_{\dot{B}{}_{q,\infty}^{s_{2}}}+\ensuremath{\sup_{0<t<T_{1}}\|N_{2}\|_{\dot{B}{}_{p,\infty}^{s_{1}}}}\\
	& \,+\left(\sup_{0<t<T_{1}}t^{\frac{\alpha}{\theta}\rho}\left\Vert E_{\alpha}(-t^{\alpha}((-\Delta)^{\theta/2}-\gamma))v_{0}\right\Vert _{\dot{B}{}_{q,\infty}^{s_{2}+\rho}}+\sup_{0<t<T_{1}}t^{\frac{\alpha}{\theta}\rho}\left\Vert E_{\alpha}(-t^{\alpha}(-\Delta^{\theta/2})\eta_{0}\right\Vert _{\dot{B}{}_{p,\infty}^{s_{1}+\rho}}\right).
\end{align*}
From hypotheses, it holds $E_{\alpha}(-t^{\alpha}((-\Delta)^{\theta/2}))\eta_{0},\eta^{1},\eta^{2}\rightarrow\eta_{0}$
and $E_{\alpha}(-t^{\alpha}((-\Delta)^{\theta/2}-\gamma))v_{0},v^{1},v^{2}\rightarrow v_{0}$
as $t\rightarrow0^{+},$ which implies that 
\begin{equation}
	\lim_{t\rightarrow0^{+}}\|N_{2}\|_{\dot{B}{}_{p,\infty}^{s_{1}}}=\lim_{t\rightarrow0^{+}}\|V_{1}\|_{\dot{B}{}_{q,\infty}^{s_{2}}}=0.\label{eq:Limitecero V1 N2}
\end{equation}
Now we prove that 
\begin{equation}
	\limsup_{t\rightarrow0^{+}}t^{\frac{\alpha}{\theta}\rho}\left\Vert E_{\alpha}(-t^{\alpha}((-\Delta)^{\theta/2}))\eta_{0}\right\Vert _{\dot{B}{}_{p,\infty}^{s_{1}+\rho}}=0.\label{eq:Limite superior eta}
\end{equation}
In fact, let $\eta_{0k}=E_{\alpha}(-\left(\frac{t}{k}\right)^{\alpha}(-\Delta)^{\theta/2})\eta_{0}$
for all $k\in\mathbb{N}.$ It follows from Lemma $\ref{Lem: Heat-Kernel-estimate en Besov}$
that $\eta_{0k}\in\dot{B}{}_{p,\infty}^{s_{1}+\rho},$ moreover, from
the hypothesis on $\eta_{0}$, we have that $\eta_{0k}\rightarrow\eta_{0}$
in $\dot{B}{}_{p,\infty}^{s_{1}}$ as $k\rightarrow\infty$. Then,
\[
\begin{aligned} & \limsup_{t\rightarrow0^{+}}t^{\frac{\alpha}{\theta}\rho}\left\Vert E_{\alpha}(-t^{\alpha}(-\Delta)^{\theta/2})\eta_{0}\right\Vert _{\dot{B}{}_{p,\infty}^{s_{1}+\rho}}\\
	& \leq\limsup_{t\rightarrow0^{+}}t^{\frac{\alpha}{\theta}\rho}\left\Vert E_{\alpha}(-t^{\alpha}(-\Delta)^{\theta/2})\left(\eta_{0}-\eta_{0k}\right)\right\Vert _{\dot{B}{}_{p,\infty}^{s_{1}+\rho}}+\limsup_{t\rightarrow0^{+}}t^{\frac{\alpha}{\theta}\rho}\left\Vert E_{\alpha}(-t^{\alpha}(-\Delta)^{\theta/2})\eta_{0k}\right\Vert _{\dot{B}{}_{p,\infty}^{s_{1}+\rho}}\\
	& \leq C\left\Vert \eta_{0}-\eta_{0k}\right\Vert _{\dot{B}{}_{p,\infty}^{s_{1}}}+C\left\Vert \eta_{0k}\right\Vert _{\dot{B}{}_{p,\infty}^{s_{1}+\rho}}\limsup_{t\rightarrow0^{+}}t^{\frac{\alpha}{\theta}\rho}\\
	& \leq C\left\Vert \eta_{0}-\eta_{0k}\right\Vert _{\dot{B}{}_{p,\infty}^{s_{1}}}\rightarrow0,\text{ as }k\rightarrow\infty.
\end{aligned}
\]
A similar argument is used to show that 
\begin{equation}
	\limsup_{t\rightarrow0^{+}}t^{\frac{\alpha}{\theta}\rho}\left\Vert E_{\alpha}(-t^{\alpha}((-\Delta)^{\theta/2}-\gamma))v_{0}\right\Vert _{\dot{B}{}_{q,\infty}^{s_{2}+\rho}}=0.\label{eq:limite cero v}
\end{equation}
Now, using estimates $(\ref{eq:Limitecero V1 N2})$, $(\ref{eq:Limite superior eta})$
and $(\ref{eq:limite cero v})$, we can choose $T_{1}>0$ such that
$CZ(T_{1})<1$ and then $N(t)=0$ for all $t\in[0,T_{1})$. From $(\ref{eq:Relaci=00003D00003D0000F3n V y N})$
we also have that $V(t)=0$ for all $t\in[0,T_{1}).$ In order to
finish the proof, we will show that $T_{1}\in(0,T]$ can be arbitrary.
For that, define
\[
T_{*}=\sup\left\{ \tilde{T};0<\tilde{T}<T,\,\eta^{1}(t)=\eta^{2}(t)\text{ in }\dot{B}{}_{p,\infty}^{s_{1}}\text{ for all }t\in[0,\tilde{T})\right\} .
\]
If $T_{*}=T$ we finish. If not, we have that $\eta^{1}(t)=\eta^{2}(t)$
for $t\in[0,T_{*})$ which implies that $\eta^{1}(T_{*})=\eta^{2}(T_{*})$
because of time continuity of $\eta^{1},\eta^{2}.$ In this case we
also have that $v^{1}(t)=v^{2}(t)$ for $t\in[0,T_{*})$ which implies
that $v^{1}(T_{*})=v^{2}(T_{*})$ because of time continuity of $v^{1},v^{2}.$
It follows from the first part of the proof, starting at $T_{*}$,
that there exists $\sigma>0$ small enough such that $\eta^{1}(t)=\eta^{2}(t)$
for $t\in[T_{*},T_{*}+\sigma),$ therefore $\eta^{1}(t)=\eta^{2}(t)$
for $t\in[0,T_{*}+\sigma),$ which contradicts the definition of $T_{*}.$

\section*{Acknowledgment}
The authors were supported by Vicerrectoría de Investigación y Extensión, UIS, CO.


\begin{thebibliography}{99}

	\bibitem{Corrias}  Corrias L.,  Perthame B., Asymptotic decay for the solutions of the parabolic–parabolic Keller–Segel chemotaxis system in
	critical spaces, \emph{Math. Comput. Modelling} \textbf{47}, 755–764 (2008)
	
	\bibitem{Duarte} Duarte-Rodríguez A.,  Ferreira L.C.F.,  Villamizar-Roa E.J., Global existence for an attraction-repulsion chemotaxis fluid model with logistic source, \emph{Discrete Contin. Dyn. Syst. Ser. B} \textbf{24}, no. 2, 423–447 (2019).
	
	\bibitem{Ferreira} Ferreira L.C.F.,  Precioso J., Existence and asymptotic behaviour for the parabolic–parabolic Keller–Segel system with
	singular data, \emph{Nonlinearity} \textbf{24}, 1433–1449 (2011).
	
	\bibitem{Hillen} Hillen T.,  Potapov A., The one-dimensional chemotaxis model: global existence and asymptotic profile, \emph{Math. Methods Appl. Sci.} \textbf{27}, no. 15, 1783–1801 (2004).
	
	\bibitem{Nagai} Nagai T., Global existence and blowup of solutions to a chemotaxis system, \emph{Nonlinear Anal.} \textbf{47}, 777–787 (2001).
	
	\bibitem{Winkler1}  Winkler M.,  Aggregation vs. global diffusive behavior in the higher-dimensional Keller–Segel model, \emph{J. Differential Equations} \textbf{248},  2889–2905 (2010).
	
	\bibitem{Winkler2}  Winkler M.,  Finite-time blow-up in the higher-dimensional parabolic-parabolic Keller-Segel system, \emph{J. Math. Pures Appl.} (9) \textbf{100}, no. 5, 748–767 (2013).
	
	\bibitem{Escudero}  Escudero C., The fractional Keller-Segel model, \emph{Nonlinearity} \textbf{19}, no. 12, 2909–2918 (2006).
	
	\bibitem{Klafter}  Klafter J.,  White B. S.,  Levandowsky M., Microzooplankton feeding behavior and the L\'evy walks in Biological Motion, W. Alt and G. Hoffmann, eds., Lecture Notes in Biomathematics, Vol. 89, Springer, Berlin, 1990. 
	
	\bibitem{Estrada}  Estrada-Rodr\'iguez G.,   Gimperlein H.,  Painter K., Fractional Patlak-Keller-Segel equations for chemotactic superdiffusion, \emph{SIAM J. Appl. Math.} \textbf{78}, no. 2, 1155-1173 (2018).
	
	\bibitem{Salem}  Salem S., Propagation of chaos for fractional Keller Segel equations in diffusion dominated and fair competition cases, \emph{J. Math. Pures Appl.} \textbf{132}, 79–132 (2019).
	
	\bibitem{Estrada2} Estrada-Rodr\'iguez G.,  Gimperlein H.,  Painter K.,  Stocek J., Space-time fractional diffusion in cell movement models with delay, \emph{Math. Models Methods Appl. Sci.} \textbf{29}, no. 1, 65-88 (2019).
	
	\bibitem{Biler}  Biler P.,  Wu G.,   Two-dimensional chemotaxis models with fractional diffusion, \emph{Math. Methods Appl. Sci.} \textbf{32}, no. 1, 112–126 (2009).  
	
	\bibitem{Zhai} Zhai Z., Global well-posedness for nonlocal fractional Keller-Segel systems in critical Besov spaces, \emph{Nonlinear Anal.} \textbf{72}, no. 6, 3173–3189 (2010).
	
	\bibitem{Bournaveas}  Bournaveas N.,   Calvez V.,  The one-dimensional Keller-Segel model with fractional diffusion of cells. \emph{Nonlinearity} \textbf{23}, no. 4, 923–935 (2010).	
	
	
	\bibitem{Shi} Shi B.,  Wang W.,  Suppression of blow up by mixing in generalized Keller-Segel system with fractional dissipation, \emph{Commun. Math. Sci.} \textbf{18}, no. 5, 1413–1440 (2020).
	
	\bibitem{Cuevas1} Azevedo J.,  Cuevas C.,  Henriquez E., Existence and asymptotic behaviour for the time-fractional Keller-Segel model for chemotaxis, \emph{Math. Nachr}. \textbf{292}, no. 3, 462–480 (2019).
	
	\bibitem{Cuevas2} Azevedo J.,  Cuevas C.,  Soto H., On the time-fractional Keller-Segel model for chemotaxis,  \emph{Math. Methods Appl. Sci.} \textbf{43}, no. 2, 769–798 (2020). 
	
	\bibitem{Planas}Carvalho-Neto P. M.,  Planas G., Mild solutions to the time fractional Navier-Stokes equations in $\mathbb{R}^n,$ \emph{J. Differential Equations} \textbf{259}, no. 7, 2948–2980 (2015).
	
	\bibitem{Jhean2} Ferreira L.C.F.,  Pérez-López J.E., Bilinear estimates and uniqueness for Navier-Stokes equations in critical Besov-type spaces, \emph{Ann. Mat. Pura Appl.} (4) \textbf{199}, no. 1, 379–400 (2020).
	
		
	\bibitem{Bony} Bony J.-M., Calcul symbolique et propagation des singularit\'es pour les \'equations aux d\'eriv\'ees
	partielles non lin\'eaires, \emph{Ann. Sci. \'Ecole Normale. Sup.}, \textbf{14}, 209–246 (1981).  
	
	\bibitem{Jhean}  Ferreira L.C.F.,  Pérez-L\'opez J.E.,  Villamizar-Roa E.J., On the product in Besov-Lorentz-Morrey spaces and existence of solutions for the stationary Boussinesq equations, \emph{Commun. Pure Appl. Anal.} \textbf{17}, no. 6, 2423–2439 (2018).
	
	\bibitem{Berg} Bergh J.,  L\"ofstr\"om J., Interpolation Spaces, Springer-Verlag, Berlin, 1976.
	
	
	
	
	
	 
	
  
	
 
	

  
  
	
	
	\bibitem{Yamazaki}  Yamazaki M., The Navier–Stokes equations in the weak-Ln space with time- dependent external force, \emph{Math. Ann.} \textbf{317} (4), 635–675 (2000).
	
	\bibitem{Meyer}  Meyer Y., Wavelets, paraproducts and Navier–Stokes equations, in: Current De- velopments in Mathematics 1996, International Press, Cambridge, MA 02238-2872, 105–212 (1999).
	

	

	
	
\end{thebibliography}
\end{document}